\documentclass[11pt, reqno ]{amsart}
\usepackage{geometry}                
\usepackage{graphicx}
\usepackage{amssymb}
\usepackage{epstopdf}
\usepackage{amsmath}
\usepackage{mathrsfs}
  \usepackage{array}
  \usepackage{amsthm}
  \usepackage{enumitem}
\usepackage[all]{xy}
\usepackage{color}
\usepackage{tikz}

\usepackage[colorlinks=true, pdfstartview=FitV,linkcolor=blue,citecolor=blue,urlcolor=blue]{hyperref}

    \advance \topmargin by -\headheight
    \advance \topmargin by -\headsep

    \evensidemargin \oddsidemargin
\newcommand{\subscript}[2]{$#1_#2$}

\numberwithin{equation}{section}

\begin{document}

 \newtheorem{theorem}{Theorem}[section]
  \newtheorem{prop}[theorem]{Proposition}
  \newtheorem{cor}[theorem]{Corollary}
  \newtheorem{lemma}[theorem]{Lemma}
  \newtheorem{defn}[theorem]{Definition}
  \newtheorem{ex}[theorem]{Example}
   \newtheorem{conj}[theorem]{Conjecture}

 \newcommand{\cx}{{\bf C}}
\newcommand{\la}{\langle}
\newcommand{\ra}{\rangle}
\newcommand{\res}{{\rm Res}}
\newcommand{\expp}{{\rm exp}}
\newcommand{\Sp}{{\rm Sp}}
\newcommand{\lgp}{\widehat{\rm G}( { F} ((t)) )}
\newcommand{\svee}{\scriptsize \vee}
\newcommand{\deff}{\stackrel{\rm def}{=}}
\newcommand{\cal}{\mathcal}

\title{Simple Vertex  Algebras Arising From Congruence Subgroups}

\author[]{Xuanzhong Dai, Bailin Song}

\address{Research Institute for Mathematical Sciences, Kyoto University, Kyoto 606-8502
JAPAN}
\email{xzdai@kurims.kyoto-u.ac.jp}

\address{School of Mathematical Sciences, University of Science and Technology of China, Hefei, Anhui 230026, P. R. CHINA}
\email{bailinso@ustc.edu.cn}

\maketitle

\maketitle

\begin{abstract}
For any congruence subgroup $\Gamma$, we consider the vertex algebra of $\Gamma$-invariant global sections of chiral de Rham complex on the upper half plane that are meromorphic at the cusps. 
We give a description of the linear structure of the $\Gamma$-invariant vertex algebra by exhibiting a linear basis determined by meromorphic modular forms, and generalize the Rankin-Cohen bracket of modular forms to meromorphic modular forms.
 We also show that the $\Gamma$-invariant vertex algebra is simple. 
\end{abstract}

\section { Introduction }\label{Section1}

The chiral de Rham complex introduced by Malikov, Schectman and Vaintrob \cite{MSV} in 1998,
 is a sheaf of vertex algebras $\Omega_X^{ch}$ defined on any smooth manifold $X$ in either the algebraic, complex analytic, or $C^\infty$ categories.
 The chiral de Rham complex has a $\mathbb Z_{\geq 0}$-grading by conformal weight, and the weight zero component coincides with the usual de Rham sheaf. Its cohomology is related to the infinite-volume limit of half-twisted sigma model.  
 
 When $X$ is the projective space $\mathbb P^n$, the vertex algebra of global sections of the chiral de Rham complex $\Gamma(X,\Omega_X^{ch})$ is computed as a module over affine Lie algebra $\widehat{\mathfrak{sl}}_{n+1}$\cite{MS}. 
 Recently the case for compact Ricci-flat K\"{a}hler manifold are studied, and the global section $\Gamma(X,\Omega_X^{ch})$ is identified as invariant elements of $\beta\gamma-bc$ system under the action of certain Lie algebra of Cartan type \cite{S2}. More specifically when $X$ is $K3$ surface, the global section is isomorphic to the simple $N=4$ algebra with central charge $6$ \cite{S1,S3}.
Let G be a connected compact semisimple Lie group and $K$ be a closed subgroup of $G$. It is natural to consider the case $X=\Gamma \backslash G /K$, where $G/K $ is a Hermitian symmetric space and $\Gamma$ is an arithmetic subgroup of $G$. 
Especially the character formula of the $\Gamma$-invariant vertex algebra of $SL(2,\mathbb R)/SO(2)$  is shown to be determined by modular forms of even weight \cite{D1,D2}, when $\Gamma$ is chosen to be a congruence subgroup of $SL(2,\mathbb Z)$. The approach is to study the global sections of the chiral de Rham complex on the upper half plane $\mathbb H$ that are holomorphic at the cusps.
 In this paper, we continue to study the chiral de Rham complex on $\mathbb H$ with a relaxed cuspidal condition, which motivates a generalization of the Rankin-Cohen bracket of modular forms as an application.

As the fractional linear transformation on $\mathbb H$ induces an $SL(2,\mathbb R)$-action on the chiral de Rham complex, we consider the vertex algebra of $\Gamma$-invariant global sections that are meromorphic at the cusps, denoted by $\mathcal M(\mathbb H,\Gamma)$.  
One of the main results is the following.
\begin{theorem} \label{main}
For any congruence subgroup $\Gamma \subset SL(2,\mathbb Z)$, the vertex algebra $\mathcal M(\mathbb H,\Gamma)$ is a simple topological vertex algebra. 
\end{theorem}

\noindent Note that $\Omega^{ch}(\mathbb H,\Gamma)$ in \cite{D2}, the vertex subalgebra of $\mathcal M(\mathbb H,\Gamma)$, corresponding to $\Gamma$-invariant sections that are holomorphic at the cusps is inherently non-simple, since the elements corresponding to modular forms of positive weights form a nontrivial ideal.

 For any partition $\lambda=(\lambda_1,\cdots,\lambda_d)$, we define $p(\lambda)=d$, and for any symbol $X=a,b,\psi$, we denote by  $X_{-\lambda}$ the expression $X_{-\lambda_1}\cdots X_{-\lambda_d}$, 
and let $\phi_{-\lambda}:=\phi_{-\lambda_1+1}\cdots \phi_{-\lambda_{d}+1}$. 
We call $(\lambda,\mu,\nu,\chi)$ a four-tuple of partitions, if $\lambda,\chi$ are partitions and $\mu,\nu$ are partitions with distinct parts. 
For any four-tuple $w$, we define the part of the four-tuple $p(w):=-p(\lambda)+p(\mu)-p(\nu)+p(\chi)$. 
The space of global sections is a free $\mathcal O(\mathbb H)$-module with a basis $a_{-\lambda}\phi_{-\mu} \psi_{-\nu}b_{-\chi}1$, when $(\lambda,\mu,\nu,\chi)$ runs through all four-tuples. 
According to Wakimoto theory \cite{W, FF, F}, the vertex algebra of global sections affords an action of affine Kac-Moody algebra $\hat{\mathfrak {sl}}_2$. Let $\{E,F,H\}$ be the standard basis of $\mathfrak{sl}_2$.
Their zeroth coefficients $E_{(0)}, F_{(0)},H_{(0)}$ induce an action of $\mathfrak {sl}_2$, which coincides with the infinitesimal action of $SL(2,\mathbb R)$.
We introduce a locally nilpotent operator which acts on the space of operators on the global sections 
 determined by 
\[
D(T)=[F_{(0)},T]+ [H_{(0)},T] b_0, \;\;\;\text{ for any operator } T.
\]

Let $\mathcal M_k(\Gamma)$ be the space of meromorphic modular forms for $\Gamma$ of level $k$ that are holomorphic on $\mathbb H$. 
The following theorem gives a linear basis of $\mathcal M(\mathbb H,\Gamma)$.  

  \begin{theorem}\label{theorem1.2} 
Let $w=(\lambda,\mu,\nu,\chi)$ be a four-tuple of partitions and $f$ be a meromorphic modular form in $\mathcal M_{2p(w)}(\Gamma)$. Then
\begin{equation}\label{1.1}
L_0(w,f):=\sum_{k=0}^\infty \frac{(\pi i)^k}{6^kk!} D^{k}(a_{-\lambda}\phi_{-\mu} \psi_{-\nu} b_{-\chi})E_2^k f \in \mathcal M(\mathbb H,\Gamma),
\end{equation}

\noindent where $E_2$ is the Eisenstein series of weight $2$. The elements $L_0(w,f)$ when $w$ runs through all four-tuples and $f$ runs through a basis of $\mathcal M_{2p(w)} (\Gamma)$, form a linear basis of $\mathcal M(\mathbb H,\Gamma)$.
\end{theorem}

\

The structure of this paper is as follows. 
In Section \ref{section2}, we review the construction of chiral de Rham complex and apply the construction to the upper half plane.
We then introduce an $SL(2,\mathbb R)$-invariant filtration on the global sections, 
and explain the relations between the $\Gamma$-fixed sections on the graded algebra and the meromorphic modular forms. 
In Section \ref{section3}, we exhibit a uniform lifting formula for meromorphic modular forms and prove Theorem \ref{theorem1.2}.
In Section \ref{section4}, we will show that the vertex operator algebra $\mathcal M (  \mathbb H,\Gamma )$ is simple for arbitrary congruence subgroup $\Gamma$.
 In Section \ref{section5}, we will give an alternative lifting formula that generalizes the lifting formulas in \cite{D1,D2}, and generalize the Rankin-Cohen bracket to meromorphic modular forms.

\

\section{Chiral De Rham Complex} \label{section2}
In this chapter, we first recall the construction of chiral de Rham complex in \cite{MSV}, which is a sheaf of vertex algebras over any smooth manifold in either algebraic, complex analytic, or $C^\infty$ settings.
We then apply the chiral de Rham complex to the upper half plane case, and introduce an $SL(2,\mathbb R)$-action on the global sections.
For any congruence subgroup $\Gamma$, we consider the vertex operator algebra consisting of  $\Gamma$-invariant global sections that are holomorphic on the upper half plane with possible poles at the cusps.

\subsection{Construction of chiral de Rham complex} \label{section2.1}

For arbitrary positive integer $N$, we recall that Heisenberg vertex algebra $V_N$ is the vacuum representation of Heisenberg algebra with basis $a^i_n, b^i_n$ for $1\leq i\leq N$ and $n\in \mathbb Z$ and center $C$ with the relation that 
\begin{equation}
[a^i_n,b^j_m] =\delta_{i,j} \delta_{n+m,0} C,
\end{equation}
and Clifford vertex algebra $ \bigwedge \nolimits_N$ is the vacuum representation of the Lie superalgebra with basis consisting of odd elements $\phi^i_n,\psi^i_n$ for $1\leq i\leq N$ and $n \in \mathbb Z$ and even center $C$ with the relation that 
\begin{equation}
\{\phi^i_n,\psi^j_m\} =\delta_{i,j} \delta_{n+m,0} C.
\end{equation}
 We consider the tensor product vertex algebra
\begin{equation} \label{2.3}
\Omega_N:=V_N\otimes \bigwedge \nolimits_N,
\end{equation}
where the Virasoro element is given by $\omega=\sum _{i=1}^N b^i_{-1}a^i+\phi^i_{-1}\psi^i$ with central charge $0$. 
The Virasoro element $\omega$ together with 
the even element $J=\sum_{i=1}^N \phi^i_0 \psi^i$ and two odd elements $Q=\sum_{i=1}^N a^i_{-1} \phi^i$ and $G=\sum_{i=1}^N b^i_{-1}\psi^i$ endows the topological vertex algebra structure of $\Omega_N$ \cite{MSV}.
 Moreover there is a vertex subalgebra in $\Omega_N$ generated by $a^i,b_{-1}^i1,\phi^i$ and $\psi^i$ for $1\leq i \leq N$ denoted by $\Omega'_N$, which is again a topological vertex algebra. One may define varies filtrations on $\Omega'_N$, and we take  $N=1$ for example. The vertex algebra $\Omega'_1$  admits a grading determined by the part of the four-tuple of partitions, namely 
 \begin{equation} \label{omega'}
 \Omega'_1=\oplus_{n\in \mathbb Z}\Omega'_1[n],
 \end{equation}
  where $\Omega'_1[n]$ is the subspace of $\Omega'_1$ spanned by elements of part $n$.

For any vector $v$ of a vertex algebra, we denote by $Y(v,z)=\sum_{n\in\mathbb Z} v_{(n)} z^{-n-1}$ the field corresponding to $v$. We call the semisimple operator $J_{(0)}$ the fermionic charge operator, and its $m$-th eigenvectors elements with fermionic charge $m$.
The chiral de Rham differential is then defined to be the odd operator $d:= -Q_{(0)}$, which increases the fermionic charge by $1$ and the square of which vanishes.
Hence the space $\Omega_N$ equipped with the chiral de Rham differential $d$ is called the chiral de Rham complex on the affine space.

Let $U$ be an open subset of an $N$-dimensional complex manifold with local coordinates $b^{1},\cdots , b^{N}$, and $\mathcal O(U)$ be the space of analytic functions on $U$. 
The localization of the chiral de Rham complex on $U$ is defined to be
\[
\Omega^{ch}(U):=\Omega_N \otimes _{\mathbb C[b^1_0,\cdots,b^N_0]} \mathcal O(U),
\]
where the action of $b^i_0$ on $\mathcal O(U)$ is simply the multiplication by $b^i$.
Then $\Omega^{ch}(U)$ is a vertex algebra generated by $a^i(z),\partial b^i(z), \phi^i(z),\psi^i(z)$ and $Y(f,z)$ for  $f\in \mathcal O(U)$, where the field $Y(f,z)$ is defined by
\begin{equation} \label{2.4}
Y(f,z)=\sum_{i= 0}^\infty \frac{\partial^i}{i!}  f(b) (\sum_{n\neq 0} b_nz^{-n})^i.
\end{equation}
We write $f(b)_{m+1}:=f(b)_{(m)}$ for the coefficient of $z^{-m-1}$ in the field $Y(f,z)$. 
These generators satisfy the following nontrivial OPEs,
\begin{gather} \label{2.5'}
a^i(z)\partial b^j(w)\sim \frac{\delta_{ij} }{(z-w)^2},\; \psi^i(z)\phi^j(w)\sim \frac{\delta_{ij}}{z-w},\\
a^i(z)f(w)\sim \frac{\frac{\partial}{\partial \, b^i} f(w)}{z-w}. \label{2.6'}
\end{gather}

\noindent If we have another coordinates $\tilde b^1,\cdots,\tilde b^N$ on $U$,
the coordinate transition equations for the generators are
\begin{align*}
\tilde a^i =& a^j_{-1} \frac{\partial b^j}{\partial \tilde b^i}  +  \frac{\partial^2 b^j}{\partial \tilde b^i \partial \tilde b^m}  \frac{\partial \tilde b^m }{\partial b^r} \phi^r \psi^j,\\
\tilde b^i_{-1}1= & \frac{\partial \tilde b^i}{\partial b^j} b^j_{-1}1,\\
\tilde \phi^i\;= &\frac{\partial \tilde b^i}{\partial b^j} \phi^i ,\\
\tilde \psi^i=& \frac{\partial b^j}{\partial \tilde b^i} \psi^j,
\end{align*}
where we use Einstein summation convention.

\subsection {$\Gamma$-invariant sections of chiral de Rham complex on the upper half plane and cuspidal conditions}

Throughout the context we will focus on the upper half plane $\mathbb H$
\[
\mathbb H:=\{ \tau \in\mathbb C \; | \text{ im }\tau >0\},
\]
and consider the vertex algebra of the global sections of chiral de Rham complex on the upper half plane
\begin{equation} \label{2.5}
\Omega^{ch}(\mathbb H) := \Omega_1 \otimes _{\mathbb C[b_0]} \mathcal O(\mathbb H).
\end{equation}
We omit the upper index and simply write $a=a^1,b=b^1,\phi=\phi^1, \psi=\psi^1$. 

According to \cite{D1} and \cite{D2}, the fractional linear transformation on $\mathbb H$
induces a right  $SL(2,\mathbb R)$-action on $\Omega^{ch}(\mathbb H)$ as vertex algebra automorphism via 
\begin{align} 
\nonumber \pi(g) a&=a_{-1}(\gamma b+\delta)^2 +2 \gamma (\gamma b+\delta) \phi _0 \psi,\\
\nonumber\pi(g) b_{-1}1&=b_{-1} (\gamma b+\delta)^{-2},\\
\label{2.6}\pi(g) \psi &=\psi_{-1} (\gamma b+\delta)^{2},\\
\nonumber\pi(g) \phi &=\phi_0 (\gamma b+\delta)^{-2},\\
\nonumber\pi(g)f(b)&=f(gb)=f\left(\dfrac{\alpha b+\beta}{\gamma b +\delta} \right),
\end{align}
for arbitrary \begin{equation} \label{2.7}
g=\begin{pmatrix}
 \alpha  & \beta\\
 \gamma &\delta
 \end{pmatrix}  \in SL(2,\mathbb R).
\end{equation}
Considering the infinitesimal action, the embedding of $\mathfrak{sl}_2$ into $\Omega^{ch}(\mathbb H)$ via
\begin{equation}\label{2.8}
E \mapsto -a,\; F \mapsto  a_{-1}b^2 +2b_0 \phi_0 \psi,\; H \mapsto -2a_{-1}b_0-2\phi_0 \psi
\end{equation}
 determines a representation of affine Kac-Moody algebra $\widehat{\mathfrak{sl}}_2$ of level $0$ \cite{W,FF,F}. Since $a_0=0$ on $\Omega'_1$ as in (\ref{omega'}), the operator $H_{(0)}$ is well defined on $\Omega'_1$. Actually it is a semisimple operator with the $-2n$-th eigenspace isomorphic to $\Omega'_1[n]$.

 For any congruence subgroup $\Gamma \subset SL(2,\mathbb Z)$, we now introduce the cuspidal conditions on the $\Gamma$-invariant global sections. Let $v$ be an arbitrary $\Gamma$-invariant global section.
 We can choose $\rho\in SL(2,\mathbb Z)$ such that $\rho (c)=\infty$.  The element
\[
\pi(\rho)v= \sum_{w=(\lambda,\mu,\nu,\chi)} a_{-\lambda} \phi_{-\mu} \psi_{-\nu} b_{-\chi} f_{w}
\]
 is invariant under $\rho^{-1} \Gamma \rho$, as the group action is a right action. 
Since $\rho^{-1} \Gamma \rho$ contains the translation matrix $\begin{pmatrix}  1 & N \\ 0 & 1 \end{pmatrix}$  for some positive integer $N$, $\pi(\rho)v$ is fixed by $\begin{pmatrix}  1 & N \\ 0 & 1 \end{pmatrix}$, which together with the following formula
 \[
\pi \left( \begin{pmatrix}  1 & N \\ 0 & 1 \end{pmatrix}  \right) \sum_{w=(\lambda,\mu,\nu,\chi)} a_{-\lambda} \phi_{-\mu} \psi_{-\nu} b_{-\chi}f_{w}(b)=\sum_{w=(\lambda,\mu,\nu,\chi)}  a_{-\lambda} \phi_{-\mu} \psi_{-\nu} b_{-\chi}f_{w}(b+N)
\]
implies that $f_{w}(b+N)=f_{w}(b)$ for any four-tuple $w$. 
  Hence we have the Fourier expansion below
 \[
  f_{w}(b)=\sum_{m\in\mathbb Z} u_{w}(m) e^{2\pi i mb/N}.
 \] 
We call $v$ is holomorphic, resp. meromorphic  at the cusp $c$ if for arbitrary four-tuple $w$, $u_{w}(m)=0$ for $m<0$, resp. $u_{w}(m)=0$ for $m<-N$ with certain positive integer $N$.
 We denote by $\Omega^{ch}(\mathbb H,\Gamma)$, resp. $\mathcal M(\mathbb H,\Gamma)$ the space spanned by the $\Gamma$-invariant vectors in $\Omega^{ch}(\mathbb H)$ that are holomorphic, resp. meromorphic at all the cusps. 
We note that $\Omega^{ch}(\mathbb H,\Gamma)$ is a vertex operator algebra \cite{D2}, and so is $\mathcal M(\mathbb H,\Gamma)$.

\subsection{Filtration induced by a semiorder}

Recall that there is a semiorder on the collection of four-tuples of partitions determined by the part of the four-tuples, namely
\begin{equation}
(\lambda,\mu,\nu,\chi)>(\lambda',\mu',\nu',\chi')\;\;\; \text{ if } \;\;\; p(\lambda,\mu,\nu,\chi)<p(\lambda',\mu',\nu',\chi').
\end{equation} 
We define a family of free $\mathcal O(\mathbb H)$-submodules:
\[
W_{m}:= \text{Span} _\mathbb C \{a_{-\lambda} \phi_{-\mu} \psi_{-\nu} b_{-\chi}f(b) \in \Omega^{ch}(\mathbb H) | p(\lambda,\mu,\nu,\chi)\geq m\},
\]
which gives an $SL(2,\mathbb R)$-invariant decreasing filtration on $\Omega^{ch}(\mathbb H)$.
The induced $SL(2,\mathbb R)$-action on the graded algebra associated to the filtration $\{ W_n\}_{n\in \mathbb Z}$ is extremely simplified. Explicitly, for any $g\in SL(2,\mathbb R)$ as in (\ref{2.7}), $f\in \mathcal O(\mathbb H)$, and four-tuple $(\lambda,\mu,\nu,\chi)$ of part $n$, we have 
\begin{equation}\label{2.10}
\pi(g) a_{-\lambda} \phi_{-\mu}\psi_{-\nu} b_{-\chi} f(b) =a_{-\lambda}\phi_{-\mu} \psi_{-\nu} b_{-\chi} (\gamma b+\delta)^{-2n}f(gb)
\mod W_{n+1}.
\end{equation}

\noindent The induced filtration on $\mathcal M(\mathbb H,\Gamma)$ is denoted by $\{ W_n(\Gamma) \}_{n\in \mathbb Z}$. 
Assume that 
\[
v=a_{-\lambda_0} \phi_{-\mu_0} \psi_{-\nu_0}b_{-\chi_0} f(b)+\sum_{\substack{w=(\lambda,\mu,\nu,\chi) \\p(w)> n}} a_{-\lambda} \phi_{-\mu} \psi_{-\nu}b_{-\chi} f_w(b) \in \mathcal M(\mathbb H,\Gamma).
\]  
By (\ref{2.10}), the function $f(b)$ must be a meromorphic modular form that is holomorphic on $\mathbb H$. Let $\mathcal M_{2n}(\Gamma)$ be the space of meromorphic modular forms of weight $2n$ with respect to $\Gamma$ that are holomorphic on $\mathbb H$.
We have the following result and the proof will be given in Section \ref{section3}.
 
  \begin{theorem} \label{theorem2.1}
 We have the short exact sequence
\begin{equation}\label{2.11}
 0\longrightarrow   W_{n+1}(\Gamma)  \longrightarrow W_n(\Gamma) \longrightarrow   \Omega_1'[n] \otimes \mathcal M_{2n}(\Gamma)   \longrightarrow 0, 
 \end{equation}
 where $\Omega_1'[n]$ is defined in (\ref{omega'}), the second arrow is the embedding map and the third arrow is the quotient map.
 \end{theorem}

 Recall that there is a unique positive definite Hermitian form on $\Omega'_1$ satisfying that $(1,1)=1$ and  \cite[\S2.3-2.4]{S2}  
\begin{align}
\label{2.12}  (a_n)^\ast= nb_{-n},\;\; (b_n)^\ast =-\frac{1}{n} a_{-n}, \;\;\;&\text{for $n\in \mathbb Z_{\neq 0}$},\\
\label{2.13}  (\psi_n)^\ast= \phi_{-n},\;\; (\phi_{n})^\ast=\psi_{-n}, \;\;\; &\text{for $n\in \mathbb Z$}.
\end{align} 
 Moreover $\Omega'_1$ is a unitary representation of the topological vertex algebra generated by $\omega,J,Q,G$, as we have
\begin{align}
Q_{(n)}^\ast= G_{(-n+1)}, \;\;\; J_{(n)} ^\ast =J_{-n},\\
\omega_{(n)} ^\ast= \omega_{(-n+2)} -(n-1) J_{(-n+1)}.
\end{align}
Let $\mathcal A$ be the vertex algebra generated by $\omega, J, Q,G$ and the following elements
\begin{gather}
 A= b_{-1}\partial a,\;\;\;C=b_{-1}a,\;\;\;
 \tilde Q= a_{-2}\phi,\;\;\; \tilde G= b_{-1}\partial \psi. 
 \end{gather}
Let $B:=\frac{1}{2}(\partial C-A)=b_{-2}a$. Thanks to the relations that
 \begin{gather}
B^\ast_{(n)}=-\frac{1}{2} A_{(-n+4)} -C_{(-n+3)},\;\;\;
C^\ast_{(n)}=C_{(-n+2)}, \\
\tilde Q_{(n)}^\ast=(-n+1) G_{(-n+2)} +\tilde G_{(-n+3)},
\end{gather}
we have the following lemma.
\begin{lemma}\label{lemma2.2}
For arbitrary $n\in \mathbb Z$, $\Omega'_1[n]$ is a unitary representation of the topological vertex algebra $\mathcal A$.
\end{lemma}

Below we will introduce another filtration on $\Omega'_1$. For any four-tuple $w=(\lambda,\mu,\nu,\chi)$, we call $p(\lambda)+p(\mu)+p(\nu)+p(\chi)$ the length of $w$ and the length of the corresponding vector $v=a_{-\lambda}\phi_{-\mu} \psi_{-\nu} b_{-\chi} 1\in \Omega'_1$ denoted by $l(w)$ or $l(v)$. 
Obviously $l(w)$ shares the same parity with $p(w)$. There is an increasing filtration $\{ \mathcal F_n \}_{n\geq 0 }$ on $\Omega'_1 $ given by
\begin{equation}
  \mathcal F_n:=\text{Span}_{\mathbb C} \{ a_{-\lambda} \phi_{-\mu} \psi_{-\nu} b_{-\chi} 1\in \Omega'_1 \; |\; 
  l(\lambda,\mu,\nu,\chi) \leq n\},
\end{equation}
The filtration is invariant under the action of $\Omega'_1[0]$.
Let $Gr\, \Omega'_1= \oplus_{n\geq 0} Gr_n\, \Omega'_1$ be the corresponding graded algebra,
where $Gr_n\, \Omega'_1=\mathcal F_{n}/\mathcal F_{n-1}$ and $\mathcal F_{-1}:=0$. For any $x \in \Omega'_1[0]$ and $m\in \mathbb Z$, we will reserve the notion $x_{(m)}$ for the induced operators on each filtered piece $Gr_n\, \Omega'_1$.

\begin{lemma} \label{lemma2.3}
 Considering $G_{(2)}$ and $\tilde G_{(3)}$ as operators on $Gr \, \Omega'_1$, we have
\begin{gather}\label{2.21}
I_{G_{(2)}}(a_{-1}^k)=I_{\tilde G_{(3)}}  (a_{-1}^k) =\emptyset,\;\;\;
I_{ G_{(2)}}  (a_{-1}^{k} \phi_{-1})=I_{\tilde G_{(3)}} (a_{-1}^{k} \phi_{-1}) =\emptyset, \\  \label{2.22}
I_{G_{(2)}} (a_{-1}^k \phi_{-1} b_{-1}^l) = a_{-1}^k \phi_{-2} \phi_{-1} b_{-1}^{l-1}1,\;\;\;\text{ for } l  \geq 1,\\     \label{2.23}
I_{\tilde{G}_{(3)}} (a_{-1}^k \phi_{-1} b_{-1}^l )= -\frac{1}{3} a_{-1}^k \phi_{-2} \phi_{-1} b_{-1}^{l-1}1,\;\;\;\text{ for } l  \geq 1,
\end{gather}
where we will not distinguish the notation between elements in $\Omega'_1$ and  their images in the graded space $Gr\, \Omega'_1$, and for any map $f$ on $Gr\, \Omega_1'$ and any four-tuple $ (\lambda,\mu,\nu,\chi)$, the set $I_f(a_{-\lambda} \phi_{-\mu} \psi_{-\nu} b_{-\chi})$ is defined to be
\[
I_f(a_{-\lambda} \phi_{-\mu} \psi_{-\nu} b_{-\chi}):=
\{ \text{ monomial } x  \in Gr\,\Omega'_1  \;| \;  
\text{ the coefficient of }  a_{-\lambda} \phi_{-\mu} \psi_{-\nu} b_{-\chi}1 \text{ in } f(x)  
\text{ equals }  1\}.
\]
\end{lemma}

\noindent  {\it Proof:}
According to Borcherds' identity, we have
\begin{align*} 
G_{(2)}=&\sum_{i\geq 0} (i-1) b_{1-i} \psi_{i} -\sum_{i\geq 0} (i+2) \psi_{-1-i}b_{i+2}\\
=&-b_{1} \psi_{0} +\sum_{i>1}(i-1) b_{1-i} \psi_{i} -\sum_{i\geq 0} (i+2) \psi_{-1-i}b_{i+2}. 
\end{align*}
Obviously the operator $b_{1}\psi_{0}$ will lower the the length strictly, and hence the induced operator on $Gr \, \Omega'_1$ acts as $0$.  The image of $G_{(2)}+b_{1}\psi_0$ 
 contains either the creation operator $b_{1-i}$ for $i>1$ or $\psi_{-1-i}$ for $i\geq 0$. Hence $a_{-1}^k1 \notin  \text{ Im } G_{(2)} \cap  Gr_{k}\, \Omega'_1$. The remaining equalities can be proved similarly.\qed

\begin{theorem} \label{theorem2.4}
For arbitrary $n \in \mathbb Z$,
$\Omega'_1[n]$ is a simple unitary module of $\Omega'_1[0]$. Especially, $\Omega'_1[0]$ is a simple  topological vertex algebra.
\end{theorem}

\noindent {\it Proof:} 
We will show that there is a unique singular vector up to a scalar in $\Omega'_1[n]$ in the sense that it is killed by all the operators below
 \begin{equation} \label{2.25}
 \omega_{(m)},J_{(k)},  Q_{(l)}, G_{(m)}, A_{(s)},B_{(s)},C_{(m)}, \tilde Q_{(k)},\tilde G_{(s)}\;\;\; \text{ for } m\geq 2, k\geq 1,l\geq 0,s\geq 3.
 \end{equation} 
  Then according to the complete reducibility property and Lemma \ref{lemma2.2}, $\Omega'_1[n]$ is a simple unitary module for $\mathcal A$ and so is it for  $\Omega'_1[0]$.
 
Take an arbitrary singular vector $w\in \Omega'_1[n]$.
Without loss of generality, we assume that $w=v + v_1+v_2+\cdots \in \Omega'_1[n] $ is homogeneous with respect to $L_0$ and $J_0$ with strictly decreasing length $l(v_i)=l(v)-2i$ for $i\geq 1$. 
Suppose there exists $m\geq 2$ such that 
\[
v=b_{-m}^k f_k+b_{-m}^{k-1} f_{k-1}+\cdots+ f_0,\;\;\; \text{ for } k\geq 1 \text{ and } f_k\neq 0,
\]
where $f_i$ is a polynomial of variables $a_{-j},b_{-l},\phi_{-s},\psi_{-t} $ with $1\leq j\leq m$, $1\leq l <m$, $s\geq 0, t\geq 1$. 
Applying the operator $A_{(m+1)}$ to $v$, we obtain that 
\[
A_{(m+1)} v+(m+1) b_{-1} (kb_{-m}^{k-1} f_k +(k-1) b_{-m}^{k-2} f_{k-1}+\cdots + f_1 ) \in \mathcal F_{l(v)-1}.
\]
Hence $A_{(m+1)} $ maps $w$ to a nonzero element in $Gr_{l(v)} \Omega'_1$ with the conformal weight strictly decreased. Similarly we can use the operator $B_{(m+1)}$ to lower the conformal weight of $w $ if 
\[
v=a_{-m}^k f_k+a_{-m}^{k-1} f_{k-1}+\cdots+ f_0, \;\;\;   \text{ for }m\geq 2, k\geq 1  \text { and } f_k\neq 0,
\]
where $f_i$ is a polynomial of variables $a_{-j},b_{-l},\phi_{-s},\psi_{-t} $ with $1\leq j < m$, $1\leq l \leq m$, $s\geq 0, t\geq 1$, and  $f_k\neq 0$. 
Hence we can assume that there are no terms $a_{-m}$ and $b_{-m}$ for $m\geq 2$ in the expression of $v$.
Then if there exists a term $\phi_{-m}$ resp. $\psi_{-m}$ for $m\geq 2$ in the expression of $v$, we can always use the operator $G_{(m)}$ resp. $Q_{(m-1)}$ to $v$ to obtain a nonzero element with strictly lower conformal weight.  
So we can assume that $v$ is a polynomial of $a_{-1},b_{-1} ,\phi_0,\phi_{-1},\psi_{-1}$.

Notice that once either $\phi_{-1}\psi_{-1} f$ or $\phi_{-1} f$ for $f\in \mathbb C[a_{-1},b_{-1}]$ appears in $v$, we can use the operator $J_{(1)}$ to decrease the conformal weight of $w$ thanks to the relations
\begin{equation}\label{2.26}
J_{(1)} \phi_{-1} f=\phi_0 f,\;\;\;   J_{(1)} \phi_{-1} \psi_{-1}f =\phi_0 \psi_{-1}f,
\end{equation}
and the fact that the induced operator $J_{(1)}$ on the graded algebra $Gr \, \Omega'_1$  kills 
\begin{equation}\label{2.27}
\phi_0 f_1+ \phi_0\phi_{-1}f_2+\phi_0 \psi_{-1} f_3+\phi_0\phi_{-1} \psi_{-1}f_4 +\psi_{-1}f_5 +f_6,\;\;\; \text{ for } f_i\in \mathbb C[a_{-1},b_{-1}]. 
\end{equation}
Thus we only need to consider the situation that $v$ is of the form (\ref{2.27}).
Since $v\in \Omega'_1[n]$ is an eigenvector for both $J_0$ and $L_0$ with fixed length $l(v)$ and fixed part $p(v)=n$, it must be a scalar multiple of one of the following cases.

\begin{enumerate}  [label=(\subscript{C}{{\arabic*}})]
 \item $ \phi_{-1} \phi_0 a_{-1}^k b_{-1}^l 1 $, \label{C1}\\
 \item $ \phi_{-1} \phi_0 \psi_{-1} a_{-1}^k b_{-1}^l 1 +c \phi_0 a_{-1}^{k+1} b_{-1}^{l+1}1 $ with $c\in \mathbb C$, \label{C2}\\
 \item $\phi_0 a_{-1}^k b_{-1}^l 1 $,  \label{C3}\\
 \item $\phi_0 \psi_{-1}  a_{-1}^k b_{-1}^l 1$ , \label{C4}\\
 \item $a_{-1}^k b_{-1}^l 1$, \label{C5}\\
 \item $\psi_{-1} a_{-1}^k b_{-1}^l 1 $, \label{C6}
\end{enumerate}
 where $k,l$ are nonnegative integers.

 Applying the operator $Q_{(0)}$ to $w$, we will obviously obtain a nonzero element in $Gr_{l(v)} \Omega'_1$ with strictly increased fermionic charge for the cases \ref{C3} ($l\geq 1$), \ref{C4}, \ref{C5} ($l\geq 1$) and \ref{C6} , as we have the relations
\begin{align}
Q_{(0)} a_{-1}^k b_{-1}^l 1=l \phi_{-1} a_{-1}^k b_{-1}^{l-1} 1&,\;\;\;
Q_{(0)} \phi_0  a_{-1}^kb_{-1}^l 1= l \phi_{-1} \phi_0  a_{-1}^k b_{-1}^{l-1} 1, \\
Q_{(0)} \phi_0 \psi_{-1}  a_{-1}^k b_{-1}^l 1  &=-\phi_0 a_{-1}^{k+1} b_{-1}^l 1+l \phi_{-1} \phi_0 \psi_{-1} a_{-1}^k b_{-1}^{l-1}  1,\\
Q{(0)}\psi_{-1} a_{-1}^k b_{-1}^l 1&= a_{-1}^{k+1}b_{-1}^l 1+l \phi_{-1} \psi_{-1} a_{-1}^k b_{-1}^{l-1} 1.
\end{align}

\noindent Similarly applying the operator $\tilde Q_{(1)}+2Q_{(0)}$ to $w$,  we will obtain a nonzero element in $Gr_{l(v)}\Omega'_1$ with strictly increased fermionic charge for the case \ref{C2}, as we have the relation that
\[
(\tilde Q_{(1)}+2Q_{(0)})(\phi_{-1} \phi_0 \psi_{-1} a_{-1}^k b_{-1}^l 1 +c \phi_0 a_{-1}^{k+1} b_{-1}^{l+1}1)=2\phi_{-1}\phi_0 a_{-1} ^{k+1}b_{-1}^l 1.
\]

For the case \ref{C3} with $l=0$ and $k\geq 1$, namely $v= a_{-1}^k\phi$,
we deduce that $G_{(2)}w$ will not be vanished
 due to Lemma \ref{lemma2.3} and the relation that  
\[
G_{(2)}  a_{-1}^k\phi=k a_{-1}^{k-1}1,\;\;\;  G_{(2)} v_i \in  \mathcal F_{k+1-2i} \text{ for } i\geq 1.
\]
The remaining situation for the case \ref{C3} is $l=k=0$, namely $v=\phi$ and hence $v_1=v_2=\cdots=0$ and $w=v=\phi$. It is simple to show that $w=\phi$ is killed by all the operators in (\ref{2.25}).  
For the case \ref{C5} with $l=0$, namely $v=a_{-1}^k1$ and hence $v_1=v_2=\cdots=0$ and $w=v=a_{-1}^k1$. We can verify that $w=a_{-1}^k1$ is killed by all the operators in (\ref{2.25}).

Now we will discuss about the final case \ref{C1}. 
If $k>0$, we apply the operator $G_{(2)}$ to $v=\phi_{-1} \phi_0 a_{-1}^k b_{-1}^l\in \mathcal F_{k+l+2}$ and obtain that
\begin{equation} \label{2.30}
G_{(2)} \phi_{-1} \phi_0 a_{-1}^k b_{-1}^l 1=k \phi_{-1} a_{-1}^{k-1} b_{-1}^l1 \in \mathcal F_{k+l}.
\end{equation}
When $l=0$,  the image of 
$G_{(2)} w$
 in $ Gr_{k} \Omega'_1$ is nonzero by (\ref{2.21}). When $l\geq 1$,  the only monomial in $\mathcal F_{k+l}$ that contributes the term $\phi_{-1} a_{-1}^{k-1} b_{-1}^l $ under the action of  $G_{(2)}$ is $\phi_{-2}\phi_{-1} a_{-1}^{k-1} b_{-1}^{l-1} 1$ by (\ref{2.22}). More explicitly we have 
\begin{equation} \label{2.31}
G_{(2)} \phi_{-2}\phi_{-1} a_{-1}^{k-1} b_{-1}^{l-1}1=\phi_{-1} a_{-1}^{k-1} b_{-1}^l1 .
\end{equation}
Notice that if the coefficient of $ \phi_{-2}\phi_{-1} a_{-1}^{k-1} b_{-1}^{l-1}1$ in $v_1$ is not equal to $-k$, then $G_{(2)}w$ is nonzero by (\ref{2.30}) and (\ref{2.31}).  Otherwise we apply $\tilde G_{(3)}$ to w, and obtain that
\[
\tilde G_{(3)} \left(\phi_{-1} \phi_0 a_{-1}^k b_{-1}^l1  -k   \phi_{-2} \phi_{-1}a_{-1}^{k-1} b_{-1}^{l-1}1 \right)
=4k \phi_{-1} a_{-1}^{k-1} b_{-1}^l 1\neq 0
\]
According to (\ref{2.23}),  the image of
$\tilde G_{(3)} w$
 in $ Gr_{k+l} \Omega'_1$ is nonzero.
If $k=0$ and $v=\phi_{-1}\phi_0 b_{-1}^l1$, we note that  $v\in \Omega'_1[l+2] \cap \mathcal F_{l+2}$ is of conformal weight $l+1$ and fermionic charge $2$, and is killed by the operators in (\ref{2.25}).
Moreover  any vector in $\Omega'_1[l+2] \cap \mathcal F_{l+2}$ of conformal weight $l+1$ and fermionic charge $2$ must be a scalar multiple of $\phi_{-1}\phi_0 b_{-1}^l1$, and $w=v=\phi_{-1}\phi_0 b_{-1}^l1$ is killed by all the operators in (\ref{2.25}). Hence there is a unique singular vector up to a scalar in each $\Omega'_1[n]$, and is proportional to $\phi_{-1}\phi_0 b_{-1}^{n-2}1$ if $n\geq 2$, $\phi$ if $n=1$ or $a_{-1}^{-n}1$ if $n\leq 0$.
\qed

\section{Uniform Lifting Formulas} \label{section3}
In this section, we study the lifting of meromorphic modular forms to $\mathcal M(\mathbb H,\Gamma)$ under the quotient map  in (\ref{2.11}). 
We will prove Theorem \ref{theorem2.1} by exhibiting a uniform lifting formula, whose function part consists of products of meromorphic modular forms and Eisenstein series of weight 2.

We first introduce a formal completion on the space of vertex operators.
Let  $\mathfrak g$ be the Lie superalgebra spanned by the odd operators $\phi_{m},\psi_{m}$ and the even operators $a_{m},b_{m}$ and the center $C$ with the Lie bracket
\begin{equation}\label{3.1}
[a_{m},b_{n}]=  \delta_{m+n,0}C ,\;\; \; [\phi_{m},\psi_{n}]=\delta_{m+n,0}C.
\end{equation}
Let $\mathcal U$ be the quotient of the universal enveloping algebra of $\mathfrak g$ by the ideal generated by $C-1$, which can be decomposed into a direct sum of conformal weight $n$ subspace $ \mathcal U_n$.
We define an increasing filtration $\{ \mathcal U_n^k\}_{k\in \mathbb Z}$ on each $\mathcal U_n$ by
 \[
 \mathcal U_{n}^k :=\sum_{i\leq k} \mathcal U_{n-i}\mathcal U_{i}.\]
Note that the collection $\{\mathcal U_n^k\}$ satisfying the property that
 \[\cap_{k\in \mathbb Z} \mathcal U_n^k =\{0\}, \;\;\; \cup_{k\in \mathbb Z} \mathcal U_n^k=\mathcal U_n
 \]
   gives a linear topology on ${\mathcal U}_n$.
Let $\bar{\mathcal U}_n$ be the completion of $\mathcal U_n$ with respect to the topology. 
Hence $\bar{\mathcal U}_n$ has a fundamental system $\{ \bar{\mathcal U}_n^k\}_{k\in \mathbb Z}$ of neighborhoods of $0$, and the direct sum $\bar {\mathcal U}:= \oplus _{n\in \mathbb Z }   \bar {\mathcal U}_n$ is a complete topological ring (see similar constructions in \cite{FZ}).
Note that $\bar{\mathcal U} $ acts on $\Omega^{ch}(\mathbb H)$.
Let $K$ be the left ideal in $\bar{\mathcal U}$ generated by elements $a_{n},b_{n},\phi_{n}, \psi_{m}$ for $n\geq 1$ and $m\geq 0$. 
Hence $\bar{\mathcal U}/K$ has the following basis
\begin{equation}\label{3.2}
a_{-\lambda}\phi_{-\mu} \psi_{-\nu} b_{-\chi}a_0^kb_0^l \;\; \; k,l\in \mathbb Z_{\geq 0}.
\end{equation}

  \noindent The operators $E_{(0)}$ and $H_{(0)}$ are well defined on $\bar{\mathcal U}/K$ and are computed below 
\begin{gather}\label{3.3}
[E_{(0)}, a_{-\lambda}\phi_{-\mu}\psi_{-\nu} b_{-\chi} a_0^k b_0^l ] =-l a_{-\lambda}\phi_{-\mu}\psi_{-\nu} b_{-\chi} a_0^k b_0^{l-1}, \\
\label{3.4}
[H_{(0)}, a_{-\lambda}\phi_{-\mu}\psi_{-\nu} b_{-\chi} a_0^k b_0^l ] =2(p(\lambda)-p(\mu)+p(\nu)-p(\chi)+k-l) a_{-\lambda}\phi_{-\mu}\psi_{-\nu} b_{-\chi} a_0^k b_0^l.
\end{gather}
We also recall the operator $D$ on $\bar{\mathcal U}/K$ defined by 
\begin{equation} \label{3.5}
D(w)=[F_{(0)},w]+[H_{(0)},w]b_0,\;\;\; \text{ for any } w\in \bar{\mathcal U}/K.
\end{equation}
It is clear that the operator $D$ will lower $H_{(0)}$-weight by $2$ and hence is locally nilpotent on $\bar{\mathcal U}/K$.

 \begin{prop} \label{prop3.1}
Assume that $w\in \bar{\mathcal U}/K$ satisfying that 
\[
[H_{(0)},w]=n w,\;\;\;[E_{(0)}, w]=0.\]
 For any $g\in SL(2,\mathbb R)$ as in (\ref{2.7}), we have 
\begin{equation} \label{3.6}
\pi(g) D^k(w) \pi(g)^{-1}=\sum_{i=0}^\infty \frac{(-1)^i}{i!} D^{k+i}(w) \gamma ^i (\gamma b+\delta)^{-2k+n-i}.
  \end{equation}
  \end{prop}
  
  \noindent {\it Proof:}
  We first show that (\ref{3.6}) holds for 
  \[
  g=\exp X\;\;\; \text{ with }
  X= 
 \begin{pmatrix}
 0 &  t \\ 0 & 0
 \end{pmatrix} \text{ and }
 \begin{pmatrix} 
0 & 0 \\  t &0
\end{pmatrix}, t\in \mathbb R.
  \]

   \noindent By definition of infinitesimal action, 
  \[
  \pi(g) D^k(w)\pi(g)^{-1} =\exp (-ad X_{(0)}) D^k(w).
   \]
 We can easily verify that (\ref{3.6}) holds for $X=\begin{pmatrix}
 0 &  1 \\ 0 & 0
 \end{pmatrix}$ and $g=\begin{pmatrix}1 &  t \\ 0 & 1 \end{pmatrix}$, as $E_{(0)}$ kills $D^k(w)$.  
For the remaining case, we have 
  \begin{align*}
  [F_{(0)}, D^k(w)]=&D(D^k(w))-[H_{(0)},D^k(w)]b_0\\
  =& D^{k+1}(w) -(n-2k) D^k(w)b_0,
\end{align*}
and 
\[
[F_{(0)}, b_0^i]= ib_0^{i+1} \mod K.
\]
  Comparing the coefficient of $t$ on both sides of (\ref{3.6}), it suffices to show
  \begin{equation} \label{3.7}
  (ad F_{(0)})^m D^k(w) =(-1)^m m! \sum_{i=0}^m \frac{(-1)^{i}}{i!} {-2k+n-i\choose m-i} D^{k+i} (w) b^{m-i}_0.
  \end{equation}
  We will prove (\ref{3.7}) by induction. Assume that (\ref{3.7}) holds for certain integer $m\geq 0$. Then 
  \begin{align*}
  (ad\, F_{(0)})^{m+1} D^k(w) =& (-1)^m m!\sum_{i=0}^m \frac{(-1)^{i}}{i!} {-2k+n-i\choose m-i}\cdot  \\&(D^{k+i+1}(w)b^{m-i}_0  -(n-2k-i-m) D^{k+i} (w) b^{m-i+1}_0 )\\
  =&(-1)^m m! \sum_{i=0}^{m+1}\frac{(-1)^{i+1}}{i!} {-2k+n-i \choose m-i+1} (m+1) D^{k+i}(w) b^{m+1-i}_0\\
  =& (-1)^{m+1} (m+1)!\sum_{i=0}^{m+1} \frac{(-1)^i}{i!} {-2k+n-i \choose m-i+1} D^{k+i}(w) b^{m+1-i}_0.
  \end{align*}

  Since the upper and lower triangular matrices with ones on the diagonal generate the whole group $SL(2,\mathbb R)$,
 it remains to show that
the right hand side of (\ref{3.6}) denoted by $F(g) D^k(w)$ satisfies the associated group law  
   \begin{equation} \label{3.8}
  F(g_2)F(g_1) D^k(w) =F(g_1g_2)D^k(w),
  \end{equation}
  where  $g_i:= \begin{pmatrix} \alpha_i & \beta_i \\ \gamma_i & \delta_i \end{pmatrix} \in SL(2,\mathbb R) $, for $i=1,2$ and $g_1g_2:=\begin{pmatrix} \alpha & \beta \\ \gamma & \delta \end{pmatrix}$. The left hand side of (\ref{3.8}) equals
  \begin{align*}
   F(g_2)F(g_1) D^k(w) = & F(g_2) \sum_{i=0}^\infty \frac{(-1)^i}{i!} D^{k+i} (w) \gamma_1^i (\gamma_1 b+\delta_1)^{-2k+n-i}\\
  = & \sum_{i,j=0}^\infty \frac{(-1)^{i+j}}{i! j!} D^{k+i+j} (w) \gamma_2^j (\gamma_2 b +\delta_2)^{-2k-2i-j+2n}\gamma_1^i (\gamma_1 g_2 b+\delta_1) ^{-2k+n-i}\\
  =& \sum_{l=0}^\infty\sum_{j=0}^l \frac{(-1)^l}{(l-j)! j!} D^{k+l}(w) \gamma_2 ^j (\gamma_2 b+\delta_2)^{-l} \gamma_1^{l-j} (\gamma b+\delta)^{-2k+n-l+j}\\
  =&\sum_{l=0}^\infty \frac{(-1)^l}{l!} \sum_{j=0}^l { l\choose j } (\gamma_2 (\gamma b+\delta))^j \gamma_1^{l-j}  D^{k+l}(w) (\gamma_2 b+\delta_2)^{-l} (\gamma b+\delta)^{-2k+n-l}\\
  =&\sum_{l=0}^\infty  \frac{(-1)^l}{l!}  D^{k+l}(w) \gamma^{l} (\gamma b+\delta)^{-2k+n-l},
  \end{align*}
  which coincides with the right hand side of (\ref{3.8}).  \qed
  
  Note that (\ref{3.6}) is equivalent to the following identities
\begin{equation}\label{3.9}
  D^m(w)= \sum_{k=m}^\infty \frac{1}{(k-m)!} \pi(g) D^k(w) \pi(g)^{-1} \gamma^{k-m} (\gamma b+\delta)^{k+m-n}, \;\;\; \text{ for } m\geq 0.
  \end{equation}

   \noindent Let $w$ be any four-tuple and $[w]$ be the corresponding vector in $\bar{\mathcal U}/K$, throughout the paper the notation $D([w])$ will be abbreviated as $D(w)$ for simplicity.  

\noindent {\it Proof of Theorem \ref{theorem1.2}:} Suppose $p(w)=n$, and so $f\in \mathcal M_{2n}(\Gamma)$. For any $g\in \Gamma$,
\begin{align*}
 \pi(g) L_0(w,f)
  =& \sum_{k\geq 0} \frac{(\pi i)^k}{6^kk!} \pi(g) D^k(w) \pi(g)^{-1} E_2^k(gb) f(gb)\\
  =& \sum_{k=0}^\infty \sum_{m=0}^k \frac{(\pi i)^m}{6^m (k-m)!m! } \pi(g) D^k(w) \pi(g)^{-1}  \gamma ^{k-m} (\gamma b+\delta)^{k+m+2n} E_2^{m} f(b)  \\
  =& \sum_{m=0}^\infty \frac{(\pi i)^m}{6^mm! }D^m(w)  E_2^{m} f(b) =L_0(w,f),
  \end{align*}
  where we use transformation equation for the Eisenstein series $E_2$ in the second equality, namely
  \[
E_2(g\tau) =(\gamma \tau+\delta)^2 E_2(\tau)-\frac{6i}{\pi} \gamma(\gamma \tau +\delta),\;\;\; \text{ for any } g=\begin{pmatrix} \alpha & \beta \\ \gamma & \delta \end{pmatrix}\in SL(2,\mathbb Z),
\]
  and (\ref{3.9}) in the third equality. Hence $L_0(w,f)\in \mathcal M(\mathbb H,\Gamma) $.
  It is clear that the set of $L_0(w,f)$ is linearly independent, and it spans the whole $\mathcal M(\mathbb H,\Gamma)$.\qed

 \noindent {\it Proof of Theorem \ref{theorem2.1}:} Using (\ref{1.1}), we can show that the quotient map (the third arrow) in (\ref{2.11}) is surjective and hence Theorem \ref{theorem2.1} holds. Indeed, for arbitrary four-tuple $w=(\lambda,\mu,\nu,\chi)$ of part $n$, and any meromorphic modular form $f\in \mathcal M_{2n}(\Gamma)$,
(\ref{1.1}) is an inverse image of $f$ with the leading term $a_{-\lambda}\phi_{-\mu}\psi_{-\nu}b_{-\chi}f(b)$ in $W_n(\Gamma)$. \qed

  \
  \

 We can define the $L_0$-operator on the tensor product space $\Omega'_1\otimes \mathcal M(\Gamma)$, where $\Omega'_1$ is defined before (\ref{omega'}) and $\mathcal M(\Gamma)=\oplus_{n\in \mathbb Z} \mathcal M_{2n} (\Gamma)$.
  For any meromorphic modular form $f\in \mathcal M_{2n}(\Gamma)$,
and arbitrary element $w=a_{-\lambda}\phi_{-\mu}\psi_{-\nu} b_{-\chi}1$ of part $n$, we  define $L_0(w,f)$ to be the right hand side of (\ref{1.1}).
 When the weight of meromorphic modular form is not twice the part of $w$, we define $L_0(w,f)$ to be $0$.
 We then extend the definition of the $L_0$ to $\Omega'_1\otimes \mathcal M(\Gamma)$ by linearity. From now on we will not distinguish four-tuple of partitions with their corresponding vectors in $\Omega'_1$.

    \begin{prop} \label{prop3.2}
  Let $w ,v$ be arbitrary two vectors in $\Omega'_1$ of part $m,n$ respectively, and let $f,h$ be arbitrary two meromorphic modular forms in $\mathcal M(\Gamma)$ of weight $2m, 2n$ respectively. We have for any $k\in \mathbb Z$
  \begin{equation} \label{3.10}
  L_0(w,f)_{(k)} L_0(v,h)- L_0\left(w_{(k)} v, fh \right)  \in  W_{m+n+1}(\Gamma).
  \end{equation}
  \end{prop}
  
  \noindent {\it Proof: } 
According to OPEs in (\ref{2.5'}) and (\ref{2.6'}), the commutator
  \[
  [a_{m},f_{(-1)}]= f'_{(m-1)}
  \] 
  will increase the part at least by $1$, and the operator $f_{(-1)}$ commutes with $X_{(m)}$ for $X\in \{ b, \psi,\phi\}$.
By Wick Theorem and (\ref{1.1}), it is easy to see that
  \begin{equation} \label{3.11}
  L_0(w,f)- w_{(-1)}f \in  W_{m+1}.
  \end{equation}
 The following identities hold modulo the space $W_{m+n+1}$
  \begin{align*}
  L_0(w,f)_{(k)} L_0(v,h)  \equiv & (w_{(-1)}f)_{(k)} v_{(-1)}h   \\
\equiv & h_{(-1)}(f_{(-1)}w)_{(k)}v \\
= & h_{(-1)} \left( \sum_{i\geq 0} f_{(-1-i)} w_{(k+i)} + \sum_{i\geq 0} w_{(-1+k-i)} f_{(i)} \right) v\\
\equiv  & h_{(-1)}f_{(-1)} w_{(k)}v  \\
\equiv & (w_{(k)}v)_{(-1)} fh,
  \end{align*}
 which implies (\ref{3.10}) due to relation (\ref{3.11}) .\qed

  \
  
  \

\section{Simplicity of vertex algebra $\mathcal M(\mathbb H,\Gamma)$}\label{section4}
 
In this section, we show that the vertex algebras $\mathcal M(\mathbb H,\Gamma)$ associated with congruence subgroups are simple. The proof is divided into two parts. The first part is an algorithm to produce the constant function $1$ from arbitrary modular form. The second part is to realize the algorithm to vertex algebra level.

  \subsection{An algorithm on modular forms} \label{section4.1}
We will introduce a linear operation on the space of modular forms $M(\Gamma)=\oplus _k M_k(\Gamma)$, which can be viewed as a special case for the modified Rankin-Cohen brackets introduced in \cite{D2}.
Especially for any $f\in  M_k(\Gamma)$, we define
\begin{equation}\label{4.1}
\{ 1, f \}:= \frac{1}{2\pi i}  \left(  f'(\tau)    - \frac{\pi i k}{6} E_2(\tau) f(\tau) \right) \in  M_{k+2}(\Gamma).
\end{equation}

\noindent As an immediate result, we have
\begin{lemma} \label{lemma4.1}
The operator $\{1, \cdot \}$ is a derivation on the ring of modular forms $M(\Gamma)$.
\end{lemma}

Recall that the discriminant function $\Delta(z)$ is the cusp form of weight $12$ for modular group $\Gamma(1)$ (eg. \cite{Z1}) with the Fourier expansion
\[
\Delta(\tau)=q\Pi_{n=1}^\infty(1-q^n) ^{24},
\]
which has the unique zero at the cusp $\infty$. 
 We also recall the Eisenstein series $E_k$ for $k> 2$ that is referred to the series 
 \[
 E_k(z)=\frac{1}{2}\sum_{\substack{c,d\in \mathbb Z \\ (c,d)=1} }  \dfrac{1}{(cz+d)^k},
 \]
 and some basic properties between $\Delta(z)$ and Eisenstein series $E_4(z)\in M_4(\Gamma(1))$ and $E_6(z)\in M_6(\Gamma(1))$.
 
 \begin{lemma}(\cite{Z1})
 We have the following relations
\begin{gather} \label{4.2}
\frac{1}{2\pi i} E_4' =\frac{1}{3} (E_2E_4-E_6),   \;\;\;  \frac{1}{ 2\pi i} E_6'=\frac{1}{2} (E_2 E_6 -E_4^2),\\
\Delta= \frac{1}{1728} (E_4^3 -E_6^2). \label{4.3}
\end{gather}
\end{lemma}

\noindent According to (\ref{4.1}) and (\ref{4.2}), we can easily compute that
\begin{equation}\label{4.4}
\{ 1, E_4\} =-\frac{1}{3} E_6,\;\;\; \{1, E_6\}=-\frac{1}{2} E_4^2.
\end{equation}
Using Lemma \ref{lemma4.1} and (\ref{4.3})-(\ref{4.4}), we have 
\begin{equation}
\{1,\Delta\}=0.
\end{equation}

 As is well-known, the ring of modular form $M(\Gamma(1))$ is generated by $E_4$ and $E_6$, namely
 \begin{equation} \label{4.6}
 M(\Gamma(1))=\mathbb C[E_4,E_6].
 \end{equation}
 Hence a modular form for $\Gamma(1)$ is either of the form
 \begin{equation} \label{4.7}
 \sum_{i=0}^k c_i \Delta^i E_4^{n-3i} \in M_{4n} (\Gamma(1)),
 \end{equation}
 or 
 \begin{equation}\label{4.8}
\sum_{i=0}^k c_i \Delta^i E_4^{n-3i} E_6 \in M_{4n+6}(\Gamma(1)).
\end{equation}
We note that the cusp $\infty$ is not a zero for both $E_4$ and $E_6$, and the term $\Delta^i E_4^{n-3i}$ vanishes at $\infty$ with multiplicity $i$. Hence the terms $\Delta^i E_4^{n-3i}$ for $0\leq i\leq k$ are linearly independent and so are the terms $\Delta^i E_4^{n-3i} E_6$ for $0\leq i \leq k$.

Below we will introduce an algorithm on arbitrary modular form by applying the operator $\{1,\cdot\}$ and a product by $E_4, E_6$ and $ \frac{1}{\Delta}$, such that the resulting function is still a modular form with the weight decreased.  

\begin{lemma}   \label{lemma4.3}
Let $f\in M_{4n}(\Gamma(1))$ be a modular form of type (\ref{4.7}). If there exists certain $j\geq 1$ such that $c_j\neq 0$, we have 
\begin{equation} \label{4.9}
0\neq \frac{1}{\Delta} \left(  E_4\{  1,f  \} +\frac{n}{3}   E_6 f   \right)      \in M_{4n-6}(\Gamma(1)).
\end{equation}
\end{lemma}

\noindent {\it Proof:} By explicit computation the resulting function in (\ref{4.9}) equals
\[
\sum_{i=0}^{k-1} (i+1) c_{i+1} \Delta^{i} E_4^{n-3i-3} E_6,
\]
which is obviously a modular form of weight $4n-6$ and nonzero if there exists $j\geq 1 $ such that $c_j\neq 0$.\qed

Now assume that $f=E_4^n$. Applying the operation $\{1,\cdot\}$ twice to $f$, we deduce that  \begin{equation} \label{4.10}
\{1, \{1,E_4^n\}\}=(-\frac{1}{6}n+\frac{1}{3} ) E_4^{n+1} +576 (n-1) \Delta E_4^{n-2} \in M_{4n+4}(\Gamma(1))
\end{equation}
due to (\ref{4.3}) and (\ref{4.4}).
When $n>2$, both the two coefficients on the right hand side of (\ref{4.10}) are nonzero. So we replace $f$ by $\{1, \{1,E_4^n\}\}$ in (\ref{4.9}), and obtain a nonzero modular form of weight $4n-2$. 
When $n=2$, (\ref{4.10}) turns into 
\[
\{1, \{1,E_4^2\}\}=576 \Delta.\]
Hence we can obtain the constant function $\frac{1}{\Delta} \{1, \{1,E_4^2\}\}=576 $. 
When $n=1$, $(\ref{4.10})$ turns into $\frac{1}{6}E_4^2$, which leads to the case for $n=2$.  According to Lemma \ref{lemma4.3} and the discussion above, we will finally arrive at the constant function $1$ by repeating the above procedures for finitely many times to arbitrary modular form of type (\ref{4.7}).

\begin{lemma}   \label{lemma4.4}
Let $f\in M_{4n+6}(\Gamma(1))$ be a modular form of type (\ref{4.8}). If there exists  certain $j\geq 1$ such that $c_j\neq 0$ we have
\begin{equation} \label{4.11}
0\neq \frac{1}{\Delta} \left(  E_4\{  1,f  \}_1 +(\frac{n}{3}+\frac{1}{2})   E_6 f   \right)      \in M_{4n}(\Gamma).
\end{equation}
\end{lemma}
\noindent {\it Proof: } Assume that $k$ is the maximal integer in the expression of $f$ such that $c_k \neq 0$. A similar calculation shows that
\begin{align*}
 \frac{1}{\Delta}  \left(  E_4\{  1,f  \} +(\frac{n}{3}+\frac{1}{2})   E_6 f   \right) 
=& \sum_{i=1}^k \left(ic_i -1728 (i-\frac{1}{2})c_{i-1}\right) \Delta^{i-1} E_4^{n-3i+3} \\
-&1728(k+\frac{1}{2})c_k \Delta^{k} E_4^{n-3k}  \neq 0
\end{align*}
which is a modular form of weight $4n$.\qed

Now assume that $f=E_4^n E_6\in M_{4n+6}(\Gamma(1))$. We have 
\[
\{1,E_4^n E_6\} =-\frac{n-1}{6} E_4^n +576 n\Delta E_4^{n-1}\neq 0,
\]
which is a modular of weight $4n+8$. When $n\geq 2$, we can apply the algorithm in Lemma \ref{lemma4.3} to obtain a nonzero modular form of weight $4n+2$. When $n=1$, 
\[
\frac{1}{\Delta} \{1,E_4 E_6\}_1 =576
\]
is already a nonzero constant function. Again we can repeat the above procedures to produce the constant $1$ from modular forms of type (\ref{4.8}).

\

\

\subsection{Operations on vertex algebra side }

In this section, we first apply the algorithm in Section \ref{section4.1} to prove the simplicity of the vertex algebra $\mathcal M(\mathbb H,\Gamma(1))$ associated to the modular group $\Gamma(1)$.  We then show that $\mathcal M(\mathbb H,\Gamma)$ is always simple for any congruence subgroup $\Gamma$.

Since the vertex algebra $\mathcal M(\mathbb H,\Gamma(1))$ is graded by the conformal weight and the fermionic charge, we start with an arbitrary $L_0$- and $J_0$-eigenvector $v$ in $\mathcal M(\mathbb H,\Gamma(1))$ denoted by
\[
v=\sum_{i=1}^m L_0(w_i,f_i )+h.o.t. ,
\]
 where the four-tuples $w_1,\cdots,w_m$ share the same part $p(w_1)$, and $h.o.t.$ refers to the liftings of meromorphic modular forms of weight greater than $2p(w_1)$.
 Assume that $w_1=(\lambda,\mu,\nu,\chi)$.  
Obviously one can use the operators
\[
L_0(b_{-1}, h_1 )_{(n)}, \; L_0(\psi_{-1}, h_2)_{(n)},\;L_0(\phi_0,h_1)_{(n)}, \;L_0(a_{-1},h_2)_{(n)},\; \text{ for }   h_1\in \mathcal M_{2}(\Gamma(1)),h_2\in \mathcal M_{-2}(\Gamma(1))
\]
to transfer $v$ into a linear combination of uniform liftings in (\ref{1.1}) with the leading term being certain meromorphic modular form of weight $0$, namely
\[
h+h.o.t,  \;\;\; \text{ for } 0\neq h \in \mathcal M_0(\Gamma(1)),
\]
where  $h.o.t$ refers to the liftings of meromorphic modular forms of positive weights.
 After multiplying $L_0(b_{-1}^6, \Delta)_{(-1)}$ for sufficiently many times to cancel the singularity of $h$ at the cusp $\infty$, the leading term turns into $b_{-1}^n f$ for certain modular form $f\in M_{2n}(\Gamma(1))$.
Hence without loss of generality, we can assume that  $v$ is an $L_0$- and $J_0$-eigenvector in $\mathcal M(\mathbb H,\Gamma(1))$ of the form
\begin{equation} \label{4.12}
v= L_0(w,f )+h.o.t. 
\end{equation}
where $f$ is a modular form of weight $2p(w)=2n>0$, and
$h.o.t.$ refers to the lifting of meromorphic modular forms of weight greater than $2n$. 
Below we will introduce vertex operations corresponding to the algorithm in Section \ref{section4.1}  which transform the leading term of $v$ to $1$. 
The multiplication of modular forms can be easily recovered by vertex operations, for example 
\[
L(u,h)_{(-1)}v=   L_0(u_{(-1)}w, hf) +h.o.t.
\] 
The only operation that is nontrivial is the operation $\{1,\cdot\}$. Below we will realize this operation to vertex algebras.

For any positive integer $n$, we denote by $w_n$ the four-tuple corresponding to $\phi_{0} \phi_{-1} \cdots \phi_{-n+1}$, namely
\[
w_n:=(\emptyset, (1,2,\cdots,n),\emptyset,\emptyset).
\]
Note that $w_n$ is the unique four-tuple of part $n$ with minimal conformal weight,
and the minimal conformal weight equals $n(n-1)/2$.
The key observation is that 
there is no four-tuple of part $n+1$ with the conformal weight $n(n-1)/2$  and fermionic charge $n$. 
 Let $v\in \mathcal M(\mathbb H,\Gamma(1))$ be an arbitrary vector as in (\ref{4.12}) with $w\in \Omega_1'[n]$ as in (\ref{omega'}).
 According to Theorem \ref{theorem2.4},
one can obtain arbitrary vector in $\Omega'_1[n]$ by successively applying vertex operators to $w$. For simplicity, we assume that there exists $u\in \Omega'_1[0]$ and $m\in \mathbb Z$ such that
\begin{equation} \label{4.13}
u_{(m)} w=\phi_{0} \phi_{-1} \cdots \phi_{-n+1}.
\end{equation}
According to Proposition \ref{prop3.2} together with the fact that there is no four-tuple of part $n+1$ with the conformal weight $n(n-1)/2$ and fermionic charge $n$, we deduce that
\begin{equation} \label{4.14}
L_0(u,1) _{(m)} L_0(w,f) - L_0(w_n,f) \in W_{n+2}(\Gamma(1)).
\end{equation}

\noindent We then apply the lifting formula (\ref{1.1}) by taking $w=w_n$ and obtain that 
\begin{equation}
L_0(w_n,f)=\phi_0\phi_{-1}\cdots \phi_{-n+1} f,
\end{equation} 
as $w_n$ is killed by the operator $D$. 

\noindent Let $u':=((1),(n+1),\emptyset,\emptyset)$ be the four-tuple corresponding to the vector $a_{-1}\phi_{-n}\in \Omega'_1[0]$. 
Applying (\ref{1.1}), the lifting of 1 with respect to $u'$ equals
\begin{equation} \label{4.17}
L_0(u',1)=
a_{-1}\phi_{-n}1 +  \frac{\pi i}{6}D(u') E_2 +\sum_{k\geq 2} \frac{(\pi i)^k}{6^k k!} D^k(u') E_2^k.
\end{equation}
To realize the operation $\{1,\cdot\}$, we will compute $L_0(u',1)_{(0)} L_0(w_n,f)$ first.
It is easy to check that
 \begin{equation} \label{4.16}
 (a_{-1} \phi_{-n}1)_{(0)} L_0(w_n,f) =(-1)^n \phi_0\phi_{-1}\cdots \phi_{-n} f'.
 \end{equation}
Applying (\ref{3.5}),
we can calculate that
\begin{align} \nonumber
D(u')=& [F_{(0)}, a_{-1}\phi_{-n}] +[H_{(0)}, a_{-1} \phi_{-n}]b_0\\
\equiv &2\sum_{i=1}^n a_{-1} \phi_{-n+i} b_{-i} +2\phi_0\phi_{-n} \psi_{-1} -2\phi_{-n-1} \mod K.  \label{4.18}
\end{align}
Hence 
by (\ref{4.16}) and (\ref{4.18}), we can compute the term $(D(u')E_2) _{(0)} L_0(w_n,f)$ as follows
\begin{align*}
(D(u')E_2) _{(0)} L_0(w_n,f) \equiv  & 2(\psi_{-1} \phi_0 \phi_{-n}E_2)_{(0)} \phi_0 \cdots \phi_{-n+1}f\\
= & 2\sum_{i= 0}^{n-1} (\phi_0\phi_{-n} E_2)_{(-1-i)}\psi_{i} \phi_0\cdots \phi_{-n+1} f   \\
= & 2\sum_{i=0}^{n-1} (-1)^ i (\phi_0 \phi_{-n} E_2) _{-i-1} \phi_0 \cdots \hat{\phi}_{-i} \cdots \phi_{-n+1}f    \\
\equiv & 2\sum_{i=0}^{n-1}\sum_{j,k\geq 0} (-1)^{i+k} {-n-1\choose k}  \phi_{-j} \phi_{-n-k}  {E_2}_{(-1-i+j+k)} \phi_0 \cdots \hat{\phi}_{-i} \cdots \phi_{-n+1} f \\
\equiv &2 \sum_{i=0}^{n-1}\sum_{j=0}^i  (-1)^{j} {-n-1\choose i-j}  \phi_{-j} \phi_{-n-i+j}   \phi_0 \cdots \hat{\phi}_{-i} \cdots \phi_{-n+1} f E_2,
\end{align*}
where the sign $``\equiv"$ means the equality holds modulo $ W_{n+2}$. 
The rightmost side of the above equality vanishes if $j\neq i$, and hence 
\begin{equation}\label{4.19}
(D(u')E_2) _{(0)} L_0(w_n,f) - 2(-1)^{n-1}n \phi_0 \phi_{-1}\cdots \phi_{-n} fE_2 \in W_{n+2} .
\end{equation}
By (\ref{4.1}), (\ref{4.17}), (\ref{4.16}) and (\ref{4.19}), we finally obtain that
\begin{align*}
L_0(u',1) _{(0)} L_0(w_n,f) 
- 2\pi i (-1)^n L_0( w_{n+1},  \{1,f\} ) \in W_{n+2}(\Gamma(1)).
\end{align*}
Therefore the operator $L_0(u',1)_{(0)} L_0(u,1)_{(-1)}$ will lead to the operation $\{1,\cdot \}$ on the modular form of the leading term, in other word we have the following
\begin{equation}
L_0(u',1)_{(0)} L_0(u,1)_{(-1)} v- 2\pi i (-1)^n L_0( w_{n+1},  \{1,f\} )  \in W_{n+2}(\Gamma(1)).
\end{equation}

\begin{theorem} \label{theorem4.5}
The vertex operator algebra $\mathcal M(\mathbb H,\Gamma(1))$ is simple.
\end{theorem}

\noindent {\it Proof: } Let $v$ be an arbitrary nonzero $L_0$- and $J_0$-eigenvector in $\mathcal M(\mathbb H,\Gamma(1))$. Thanks to the algorithm on modular forms, we can assume that the modular form of the leading term of $v$ is the constant function $1$, namely
\begin{equation} 
v=L_0(w,1) +h.o.t,
\end{equation}
where $h.o.t.$ refers to liftings of meromorphic modular forms of positive weights. According to Theorem \ref{theorem2.4}, we can assume that there exists $u\in \Omega'_1[0]$ and $m\in \mathbb Z$ such that $u_{(m)}v=1$. Hence by Proposition \ref{prop3.2} we have 
\[
L_0(u,1)_{(m)} v= 1+h.o.t.
\]
Since both the conformal weight and fermionic charge of $1$ are $0$, this forces $h.o.t.$ to be $0$. Therefore one can obtain the vacuum vector $1$ from arbitrary nonzero element via certain vertex operations, and thus $\mathcal M(\mathbb H,\Gamma(1))$ must be simple. \qed

Recall that the principal congruence subgroup of level $N$ in $\Gamma(1)$ denoted by  $\Gamma(N) $,  is the kernel of the homomorphism $\Gamma(1) \longrightarrow SL(2, \mathbb Z/N\mathbb Z)$ induced by the modulo $N$ morphism $\mathbb Z \longrightarrow \mathbb Z / N\mathbb Z$. 
Hence $\Gamma(N)$ is a normal subgroup of $\Gamma(1)$, and $\Gamma(1)$ can be decomposed as a disjoint union of right cosets,
\begin{equation} \label{4.22}
\Gamma(1) =\cup _{i=1}^k \Gamma(N) g_i, \;\;\;\text{ for } g_i \in \Gamma(1) \text{ and } g_k =\begin{pmatrix} 1 & 0 \\ 0& 1 \end{pmatrix}
\end{equation}

\begin{theorem}
For any congruence subgroup $\Gamma\subset \Gamma(1)$, the vertex operator algebra $\mathcal M(\mathbb H,\Gamma)$ is simple.
\end{theorem}

\noindent  {\it Proof: }  
Take an arbitrary nonzero $L_0$- and $J_0$-eigenvector $v\in \mathcal M(\mathbb H,\Gamma)$. 
Similar to the discussion before (\ref{4.12}), we can assume 
\begin{equation}
v=L_0(1, f) +  h.o.t. \in \mathcal M(\mathbb H,\Gamma),  
\end{equation}
where $f\in \mathcal M_0(\Gamma)$, and $L_0(1,f) =f$.
Assume that there exists $N>0$ such that $\Gamma \supset \Gamma(N)$ and
the right coset decomposition of $\Gamma(1)$ in terms of $\Gamma(N)$ as in (\ref{4.22}). 
Let 
\[
\tilde f(b)= \Pi_{i=1}^{k-1} \pi(g_i) f(b)=\Pi_{i=1}^{k-1} f(g_ib).
\]
Since $f(b) \tilde f(b)=\Pi_{i=1}^k f(g_i b)$ is contained in $\mathcal M_0(\Gamma(1))\subset \mathcal M_0(\Gamma)$, we have 
\[
(\tilde f(gb)-\tilde f(b)) f(b)=0\;\;\;\text{ for any } g\in \Gamma.
\]
As $f$ is a holomorphic function not identically zero, the holomorphic function $\tilde f(gb)-\tilde f(b)$ must be vanished on $\mathbb H $ for any $g\in \Gamma$
and hence $\tilde f(b)$ is contained in both $\mathcal M_0(\Gamma)$ and  $\mathcal M(\mathbb H,\Gamma)$.
Thus the leading term of
$\tilde f(b)_{(-1)} v$ equals 
\[\Pi_{i=1}^k f(g_i b)=L_0(1,\Pi_{i=1}^k f(g_i b)) \in \mathcal M(\mathbb H,\Gamma(1)).\] 
Thanks to the algorithm mentioned before, we finish the proof by transforming the leading term to the vacuum vector $1$. \qed

 \section{Invariant global sections and Jacobi-like forms} \label{section5}
 In this section, we will construct invariant global sections distinct from (\ref{1.1}), which can be viewed as a generalization of the lifting formulas in \cite{D2}.
 We will also explore the relations between invariant sections and Jacobi-like forms, and generalize the Rankin-Cohen brackets to meromorphic modular forms.
 
 \subsection{Alternative version of Lifting formulas}

We first introduce a simple lemma concerning the derivatives for meromorphic forms.

\begin{lemma} \label{lemma5.1}
Let $n\in \mathbb Z_{\geq 0}$, and $f\in \mathcal M_{-n}(\Gamma)$. For any $0\leq k\leq n$, we have
\begin{equation} \label{5.1}
\frac{(-1)^k (n-k)!} {k!} f^{(k)} (g\tau)=\sum_{m=0}^k \dfrac{\gamma ^{k-m} (\gamma \tau +\delta)^{-n+k+m}} {(k-m)!} \frac{(-1)^m (n-m)!}{m!} f^{(m)}(\tau),
\end{equation}
and $f^{(n+1)}\in \mathcal M_{n+2}(\Gamma)$.
  \end{lemma}
  
  The lemma shows that the modularity property for $f\in \mathcal M_{-n}(\Gamma)$ mutates, when taking $(n+1)$-th derivatives.
  We denote by $c_{-n,m}:=\frac{(-1)^m (n-m)!}{m!}$ for $0\leq m \leq n$.
  
  \begin{lemma} \label{lemma5.2}
Let $w$ be a four-tuple of part $-n$ with $n\geq 0$, and let $f\in \mathcal M_{-2n}(\Gamma)$.  For any $g\in \Gamma$
\begin{equation}
\pi(g) \sum_{m=0}^{2n} c_{-2n,m} D^m(w) f^{(m)}(b)-\sum_{m=0}^{2n} c_{-2n,m} D^{m}(w) f^{(m)}(b)\in W_{n+1}
\end{equation}
  \end{lemma}
  
  \noindent {\it Proof:} According to Proposition \ref{prop3.1} and (\ref{5.1}), we have
  \begin{align*}
  \pi(g)\sum_{m=0}^{2n} c_{-2n,m} D^m(w) f^{(m)} (b)=& \sum_{m=0}^{2n} c_{-2n,m} \pi(g) D^m(w) \pi(g)^{-1} f^{(m)}(gb)\\
  =& \sum_{m=0}^{2n} \sum_{i=0}^\infty \sum_{k=0}^m \frac{(-1)^{i+k} (2n-k)!}{i! (m-k)! k!} D^{m+i}(w) \gamma^{m+i-k}\\
  & \cdot (\gamma b+\delta)^{-m+k-i} f^{(k)}(b)\\
  =& \sum_{k=0}^{2n} \sum_{i=k} ^\infty \sum_{m=k}^{\min\{2n, i\}} \frac{(-1)^{i-m+k} (2n-k)!}{(i-m)!(m-k)! k! } D^i(w) \gamma^{i-k} \\
  & \cdot  (\gamma b+\delta)^{-i+k} f^{(k)}(b).
  \end{align*}
  
  \noindent The part corresponding to $D^i(w)$ for $0\leq i \leq 2n$ equals
  \begin{align*}
 &\sum_{k=0}^{2n} \sum_{i=k}^{2n}\sum_{m=k}^{i} \frac{(-1)^{i-m+k} (2n-k)!}{(i-m)!(m-k)! k! } D^i(w) \gamma^{i-k} (\gamma b+\delta)^{-i+k} f^{(k)}(b)\\
 =&\sum_{k=0}^{2n}\sum_{i=k}^{2n}\sum_{m=0}^{i-k} \frac{(-1)^{m}}{(i-k-m)!m!} \frac{(-1)^i(2n-k)! }{k!} D^i(w) \gamma^{i-k} (\gamma b+\delta)^{k-i} f^{(k)}(b)\\
 =& \sum_{k=0}^{2n} \frac{(-1)^k (2n-k)!}{k!} D^k(w) f^{(k)}(b).
   \end{align*}

\noindent Since the part of $D^{2n+k}(w)$ equals $n+k$ for $k\geq 1$, the remaining terms must be contained in $W_{n+1} $.\qed

Thanks to Lemma \ref{lemma5.2}, we can assume that the first $(2n+1)$-terms of the lifting of $f\in \mathcal M_{-2n}(\Gamma)$ equals
\[\sum_{m=0}^{2n} c_{-2n,m} D^{m}(w) f^{(m)}(b).\]
Notice that for $f\in \mathcal M_{-2n}(\Gamma)$, the modularity property of $f^{(k)}$ mutates when $k=2n+1$ by Lemma \ref{lemma5.1}. 
 To derive successive terms, we need to introduce a series of functions $\{h^f_m\}_{m\geq 0}$ depending on the meromorphic form $f\in \mathcal M_{-2n}(\Gamma)$, where $h^f_m$ is a holomorphic function on $\mathbb H$ defined by
  \begin{equation}\label{5.3}
  h^f_m(b):=\frac{\pi i}{6m!} \sum_{k=0}^{2n} \frac{1}{k! (2n-k+m+1)!} f^{(k)}(b) E^{(2n+m-k)}_2(b).
  \end{equation}
  
  \begin{lemma} 
  Let $n\in \mathbb Z_{\geq 0}$, and $f\in \mathcal M_{-2n}(\Gamma)$. We have 
  \begin{equation}\label{5.4}  
 h^f_k (b)'=(2n+k+2)(k+1) h^f_{k+1} (b)+\frac{\pi i}{6(2n)! k!(k+1)!} f^{(2n+1)}(b) E_2^{(k)}(b).
  \end{equation}
  \end{lemma}

As a preliminary lemma, we will show that $h^f_0$ is a suitable choice for the coefficient of $D^{2n+1}(w)$.
\begin{lemma} \label{lemma5.4}
Let $w$ be a four-tuple of part $-n$ with $n\geq 0$, and let $f\in \mathcal M_{-2n}(\Gamma)$.  For any $g\in \Gamma$, we have
\begin{align*}
\pi(g)( \sum_{m=0}^{2n} &c_{-2n,m} D^m(w) f^{(m)}(b) +D^{2n+1}(w)h_0^f(b) )\\
&-(\sum_{m=0}^{2n} c_{-2n,m} D^{m}(w) f^{(m)}(b) +D^{2n+1}(w) h_0^f(b)) \in   W_{n+2}.
\end{align*}
  \end{lemma}

\noindent {\it Proof:} 
Considering the function part corresponding to $D^{2n+1}(w)$, it suffices to show
\begin{equation}\label{5.5}
(\gamma b+\delta)^{-2n-2} h^f_0(gb) =h^f_0(b) +\sum_{m=0}^{2n} c_{-2n,m} \frac{(-1)^m}{(2n-m+1)!} \gamma^{2n-m+1} (\gamma b+\delta)^{-m-1} f^{(m)} (gb)
\end{equation}

\noindent Note that the Eisenstein series $E_2$ satisfies that
\begin{equation} \label{5.6}
\dfrac{E_2^{(k)}(g \tau)}{k! (k+1)!} =\sum_{m=0}^k  \dfrac{ \gamma^{k-m} (\gamma \tau +\delta)^{m+k+2}}{(k-m)!}\dfrac{E_2^{(m)}(\tau)}{m!(m+1)!}-\frac{6i}{\pi} \dfrac{\gamma^{k+1} }{(k+1)!} (\gamma \tau +\delta)^{k+1} \;\;\;\;\text{ for } k\geq 0,
\end{equation}
which can be easily proved by induction.
  Using (\ref{5.1}) and (\ref{5.6}), we compute that $h^f_0(gb)=L_1+L_2$, where 
  \[
  L_1= \sum_{k=0}^{2n} \sum_{m=0}^k \sum_{l=0}^{2n-k} \dfrac{\pi i(-1)^{k+m} (2n-m)! }{6(k-m)! (2n-k-l)! l! (l+1)!m! } \gamma^{2n-m-l} (\gamma b+\delta)^{m+l+2} f^{(m)}(b) E_2^{(l)}(b),
  \]
  and 
  \[
  L_2=\sum_{k=0}^{2n} \sum_{m=0}^k \dfrac{(-1)^{k+m} (2n-m)!}{(k-m)!(2n-k+1)!m!} \gamma^{2n-m+1} (\gamma b+\delta)^{m+1} f^{(m)}(b).
\]  
  After a change of variable $m+l \rightarrow s$, we have
  \begin{align}  \nonumber
  L_1=&\frac{\pi i} {6}\sum_{s=0}^{2n} \sum_{m=0}^s \sum_{k=m}^{2n-s+m} \dfrac{(-1)^{k+m} (2n-m)! }{(k-m)! (2n-k-l)! l! (l+1)!m! } \gamma^{2n-m-l} (\gamma b+\delta)^{m+l+2}\\
  &\cdot f^{(m)}(b) E_2^{(l)}(b)\\   \nonumber
  =&\frac{\pi i}{6}\sum_{s=0}^{2n} \sum_{m=0}^s \sum_{k=0}^{2n-s} \frac{(-1)^{k}}{k!(2n-s-k)!} \frac{(2n-m)!}{(s-m)!(s-m+1)! m!} \gamma^{2n-s} (\gamma b+\delta)^{s+2}\\
  &\cdot f^{(m)}(b) E^{(s-m)}_2(b)\\ \label{5.9}
  =&\frac{\pi i}{6} \sum_{m=0}^{2n} \frac{1}{(2n-m+1)!m!} (\gamma b+\delta)^{2n+2} f^{(m)}(b) E^{(2n-m)}_2(b).
  \end{align}
  Changing the summation order of $L_2$, we have
  \begin{align} \nonumber
  L_2=&\sum_{m=0}^{2n}\sum_{k=m}^{2n}\dfrac{(-1)^{k+m} (2n-m)!}{(k-m)!(2n-k+1)!m!} \gamma^{2n-m+1} (\gamma b+\delta)^{m+1} f^{(m)}(b)\\   \nonumber
  =& \sum_{m=0}^{2n}\sum_{k=0}^{2n-m} \frac{(-1)^{k} }{k! (2n+1-m-k)! } \frac{(2n-m)!}{m!}\gamma^{2n-m+1} (\gamma b+\delta)^{m+1} f^{(m)}(b)\\ \nonumber
  =& \sum_{m=0}^{2n} \frac{(-1)^{m}}{(2n-m+1)!} \frac{(2n-m)!}{m!} \gamma ^{2n-m+1} (\gamma b+\delta)^{m+1} f^{(m)}(b)\\  \label{5.10}
  =& \sum_{m=0}^{2n} \frac{(-1)^{m}}{(2n-m+1) m!} \gamma ^{2n-m+1} (\gamma b+\delta)^{m+1} f^{(m)}(b)
  \end{align}
  
  \noindent The second summand on the right hand side of (\ref{5.5}) equals
  \begin{align}\nonumber
 & \sum_{m=0}^{2n}\sum_{i=0}^m \frac{(-1)^m}{(2n-m+1)! (m-i)! } \frac{(2n-i)!}{i! (2n)!} \gamma^{2n-i+1} (\gamma b+\delta)^{-2n+i-1} f^{(i)}(b)\\ \nonumber
  = & \sum_{i=0}^{2n} \sum_{m=i}^{2n}  \frac{(-1)^m}{(2n-m+1)! (m-i)! } \frac{(2n-i)!}{i! (2n)!} \gamma^{2n-i+1} (\gamma b+\delta)^{-2n+i-1} f^{(i)}(b)\\ \label{5.11}
  =& \sum_{i=0}^{2n} \frac{(-1)^i}{(2n+1-i)! i!} \gamma^{2n-i+1} (\gamma b+\delta)^{-2n+i-1} f^{(i)}(b)
  \end{align}
  
 \noindent  Combing (\ref{5.9}), (\ref{5.10}), (\ref{5.11}) and the definition of  $h^f_0$, we will obtain (\ref{5.5}).\qed

Now we have the following lifting theorem.
 \begin{theorem}\label{theorem5.5}
  Let $w=(\lambda,\mu,\nu,\chi)$ be a four-tuple of partitions such that $p(w)=-n$ for $n\geq 0$, and $f$ be a meromorphic modular form for $\Gamma$ of weight $-2n$. Then
\begin{equation}\label{5.12}
L_1(w,f):=\sum_{m=0}^{2n} c_{-2n,m} D^{m}(w) f^{(m)} +\sum_{m=0}^\infty D^{2n+1+m}  (w) h^f_m \in \mathcal M(\mathbb H,\Gamma).
\end{equation}
  \end{theorem}
  \noindent {\it Proof of Theorem \ref{theorem5.5}:} Thanks to Lemma \ref{lemma5.2} and \ref{lemma5.4},  we only need to compare the coefficients of $D^{2n+1+s}(w)$ for $\pi(g)L_1(w,f)$ and $L_1(w,f)$, which is equivalent to the identity below
    \begin{align}\nonumber
  \sum_{m=0}^{2n} c_{-2n,m}& \frac{(-1)^{s-m+1}}{(2n+s-m+1)!} \gamma^{2n+s-m+1} (\gamma b+\delta)^{-m-s-1} f^{(m)}(gb)  \\ \label{5.13}
 & +\sum_{k=0}^s \frac{(-1)^{s-k}} {(s-k)!} \gamma^{s-k} (\gamma b+\delta)^{-2n-s-k-2} h^f_k(gb) =h^f_s(b).
  \end{align}
  \noindent Obviously the case $s=0$ for (\ref{5.13}) is exactly (\ref{5.5}).
   Assume (\ref{5.13}) holds for $s\leq l$. Let 
  \begin{equation} 
  L_1=\sum_{m=0}^{2n}  c_{-2n,m} \frac{(-1)^{l-m+1}}{(2n+l-m+1)!} \gamma^{2n+l-m+1} (\gamma b+\delta)^{-m-l-1} f^{(m)}(gb), 
  \end{equation}
  and 
  \begin{equation}
  L_2=\sum_{k=0}^s \frac{(-1)^{l-k}} {(l-k)!} \gamma^{l-k} (\gamma b+\delta)^{-2n-l-k-2} h^f_k(gb). 
  \end{equation}
We will calculate the derivative of $L_1$ and $L_2$, respectively.
  \begin{align*}
  L_1'=&(l+1)(l+2n+2) \sum_{m=0}^{2n} \frac{c_{-2n,m} (-1)^{m+l}}{(2n+l-m+2)!} \gamma^{2n+l-m+2}(\gamma b+\delta)^{-m-l-2} f^{(m)}(gb)\\
  &+\frac{(-1)^{l+1}}{(2n)! (l+1)!} \gamma^{l+1} (\gamma b+\delta)^{-l-1} f^{(2n+1)}(b),
  \end{align*}
\begin{align*}  
  L_2'=& \sum_{k=0}^l \frac{(-1)^{l-k}}{(l-k)!} \gamma^{l-k+1} (\gamma b+\delta)^{-2n-l-k-3} (2n-l-k-2) h^f_k(gb)\\
  &+\sum_{k=0}^l \frac{(-1)^{l-k}}{(l-k)!} \gamma^{l-k} (\gamma b+\delta)^{-2n-l-k-2} h^f_k(gb)'.
  \end{align*}
 \noindent We denote by $L_3$ and $L_4$ the first and the second summand of $L_2'$ respectively. 
  
    \begin{align*}
  L_4=& \sum_{k=0}^l \frac{(-1)^{l-k}}{(l-k)!} \gamma^{l-k} (\gamma b+\delta)^{-2n-l-k-4} \\
 & \left( (2n+k+2)(k+1) h^f_{k+1}(gb) 
+ \frac{\pi i}{6(2n)!k!(k+1)!} f^{(2n+1)}(gb) E_2^{(k)}(gb) \right)  \\
    =&\sum_{k=0}^l \dfrac{(-1)^{l-k} (2n+k+2)(k+1) }{(l-k)!} \gamma^{l-k} (\gamma b+\delta)^{-2n-l-k-4} h^f_{k+1}(gb) \\
  &+\frac{\pi i}{6}\sum_{k=0}^l \dfrac{(-1)^{l-k}}{(l-k)! (2n)! k! (k+1)! } \gamma^{l-k} (\gamma b+\delta)^{-2n-l-k-4} f^{(2n+1)}(gb) E_2^{(k)}(gb)
  \end{align*}
  \noindent Again we denote by $L_5$ and $L_6$ the first and second summand of $L_4$ respectively.
We have   
  \begin{align*}
  L_3+L_5=&\sum_{k=0}^{l} \frac{(-1)^{l-k}}{(l-k)!} \gamma^{l-k+1} (\gamma b+\delta)^{-2n-l-k-3} (-2n-l-k-2) h^f_k(gb)\\
  &+\sum_{k=0}^l \dfrac{(-1)^{l-k} (2n+k+2)(k+1)}{(l-k)!}\gamma^{l-k} (\gamma b+\delta)^{-2n-l-k-4} h^f_{k+1}(gb)\\
  =& \sum_{k=1}^l \left(  \frac{(-1)^{l-k+1} (2n+l+k+2) } {(l-k)!}  +\frac{(-1)^{l-k+1} (2n+k+1) k  } {(l-k+1)!}     \right)\\
  &\cdot \gamma^{l-k+1} (\gamma b+\delta)^{-2n-l-k-3} h^f_k(gb)\\
  &+ \frac{(-1)^{l+1} (2n+l+2)} {l!} \gamma^{l+1} (\gamma b+\delta)^{-2n-k-3} h^f_0(gb)\\
  &+ (2n+l+2)(l+1) (\gamma b+\delta)^{-2n-2l-4} h^f_{l+1}(gb)\\
  =&(l+1)(2n+l+2) \sum_{k=0}^{l+1} \frac{(-1)^{l-k+1}}{ (l-k+1)!} \gamma^{l-k+1} (\gamma b+\delta)^{-2n-l-k-3}  h^f_k(gb), 
  \end{align*}
  and also by Lemma \ref{lemma5.2} and (\ref{5.6})
  \begin{align*}
  L_6=& \sum_{k=0}^l \dfrac{(-1)^{l-k}}{(l-k)! (2n)!} \gamma ^{l-k} (\gamma b+\delta)^{-l-k-2 } f^{(2n+1)}(b)\\
  &\cdot \left( \frac{\pi i}{6} \sum_{m=0}^k \dfrac{\gamma^{k-m} (\gamma b+\delta)^{k+m+2}} {(k-m)!} \dfrac{E_2^{(m)}}{m!(m+1)!} +\dfrac{\gamma ^{k+1}}{(k+1)!} (\gamma b+\delta)^{k+1} \right)\\
  =&\frac{\pi i}{6} \sum_{m=0}^{l} \sum_{k=m}^l \frac{(-1)^{l-k}}{(l-k)! (k-m)!} \frac{1}{(2n)! m! (m+1)!} \gamma^{l-m }(\gamma b+\delta)^{-l+m} f^{(2n+1)}(b) E^{(m)}_2(b) \\
  &+ \sum_{k=0}^l \frac{(-1)^{l-k}}{(l-k)! (k+1)!} \frac{1}{(2n)!} \gamma^{l+1} (\gamma b+\delta)^{-l-1}\\
  =& \frac{\pi i}{6}\frac{(1)}{(2n)! l! (l+1)!}  f^{(2n+1)}(b) E_2^{(m)}(b)\\
  &+ \frac{(-1)^l}{(l+1)! (2n)!} \gamma^{l+1} (\gamma b+\delta)^{-l-1} f^{(2n+1)}(b)
  \end{align*} 
  
\noindent  Notice that the relation $L_1' +L_3+L_5+L_6=h^f_l(b)'$ together with (\ref{5.4}) leads to (\ref{5.13}) for $s=l+1$, and induction works.\qed

  Note that the formula $(1.3)$ in \cite{D2} can be viewed as a special case of (\ref{5.12}),
and the formula (1.2) in \cite{D2} also gives a lifting formula for meromorphic modular form, namely when $f\in \mathcal M_{2n}(\Gamma)$ and $p(w)=n$ with $n>0$
   \begin{equation} \label{5.16}
  L_1(w,f):= \sum_{m=0}^\infty \frac{(2n-1)!}{m! (m+2n-1)!} D^{m}(w) f^{(m)}(b) \in \mathcal M(\mathbb H,\Gamma).
  \end{equation}
  Obviously the elements $L_1(w,f)$ in (\ref{5.12}) and (\ref{5.16}) when $w$ runs through all four-tuples, and $f$ runs through a basis of meromorphic modular form of weight $2p(w)$, form a linear basis of $\mathcal M(\mathbb H,\Gamma)$.

 For arbitrary $f\in \mathcal M_{2n}(\Gamma), n\in \mathbb Z$, we define 
  \begin{equation}
  f^{[m]}:=
  \begin{cases} 
  f^{(m)}  \;\;\; & \text{ if } n>0, \text{ or } n\leq 0 \text{ and } 0\leq m\leq -2n;\\
  h^f_{m+2n-1} & \text{ if } n\leq 0, m\geq -2n+1,
  \end{cases}
  \end{equation}
  and define $c_{n, m}:= 1$ when $n<0, m \geq -n +1$ and $c_{n,m}:= \frac{(n-1)!}{m! (m+n-1)!}$ when $n>0$.
We can rewrite $L_1(w,f)$ in (\ref{5.12}) and (\ref{5.16}) as
\begin{equation}\label{5.18}
L_1(w,f)=\sum_{m=0}^\infty c_{2n,m} D^{m}(w) f^{[m]} \in \mathcal M(\mathbb H,\Gamma).
\end{equation}

  \
  
  \

\subsection{The Generalized Rankin-Cohen bracket}

The Rankin-Cohen bracket is a family of bilinear operations, which sends two modular forms $f$ of weight $k$ and $h$ of weight $l$, to a modular form $[f,h]_n$ of weight $k+l+2n$. Let $\Gamma\subset SL(2,\mathbb Z)$ be a congruence subgroup, and $f\in M_{k}(\Gamma)$, $h\in M_{l}(\Gamma)$, then the $n$-th Rankin-Cohen bracket is given by
\begin{equation}\label{5.19}
[f,h]_n=\frac{1}{(2\pi i)^n }\sum_{r+s=n} (-1)^r {n+k-1\choose s} {n+ l-1\choose r} f^{(r)}(\tau) h^{(s)}(\tau).
\end{equation}
 Therefore the graded vector space $M (\Gamma) =\oplus_{i\geq 0}M_{i}(\Gamma)$ possesses an infinite family of bilinear operations $[\; ,\; ]_n: M_\ast \otimes M_{\ast}\mapsto M_{\ast+\ast+2n}$, with the $0$th bracket the usual multiplication. 
 The modularity of the Rankin-Cohen bracket was first proved by Cohen \cite{C}, using the Cohen-Kuznetsov lifting from modular forms to Jacobi-like forms (see \cite{Z2} and \cite{CMZ}).
 Let $R$ be the ring of all holomorphic functions on $\mathbb H$ bounded by a power of $(|\tau|^2+1)/\text{Im}(\tau)$, namely
 \[
 R:=\{ f\in \mathcal O(\mathbb H) \;| \;  |  f(\tau)| \leq C (|\tau|^2+1)^l/\text{Im}(\tau)^l \;\;\text{ for some } l ,C>0\}.
 \]
 A Jacobi-like form of weight $k$ for $\Gamma$  is defined to be a formal power series $\Phi(\tau,X)\in R[[X]]$, such that  
  for any $g= \begin{pmatrix} \alpha & \beta \\ \gamma & \delta \end{pmatrix} \in \Gamma$,   
  \begin{equation}\label{5.20}
\Phi\left(\dfrac{\alpha \tau +\beta}{\gamma \tau +\delta}, \dfrac{X}{ (\gamma \tau +\delta)^2}\right) =(\gamma \tau +\delta)^k e^{\gamma  X/(\gamma \tau +\delta)} \Phi(\tau,X).
\end{equation}

We introduce a family of variables $\{X_{n,m}\}_{n\in \mathbb Z, m\in \mathbb Z_{\geq 0}}$ labeled by $\mathbb Z \times \mathbb Z_{\geq 0}$.
Let $\mathcal X_n$ be the space spanned by formal series 
$
\sum_{m=0}^\infty f_m X_{n,m}, f_m\in R.
$ We equip an $SL(2,\mathbb R)$-action on $\mathcal X_n$ via
\[f_mX_{n,m} \mapsto \sum_{i=0}^\infty \frac{(-1)^i}{i!}   \gamma^i (\gamma b+\delta)^{-2k-2n-i} f_m(gb) X_{n,m+i}, \;\;\;\text{ for } g= \begin{pmatrix} \alpha & \beta \\ \gamma & \delta \end{pmatrix} \in SL(2,\mathbb R).\]

\begin{theorem} \label{theorem5.6}
There is a one-one correspondence between $\Gamma$-fixed points in $\mathcal X_n$ and Jacobi-like forms of weight $2n$ for $\Gamma$, by sending $\sum_{m=0}^\infty f_m X_{n,m}$ to $\sum_{m=0}^\infty f_m X^m$.
\end{theorem}

\noindent {\it Proof: } Let $L_n=\sum_{m=0}^\infty  f_m X_{n,m} $. Then we have
\begin{align*}
\pi(g)L_n =&\sum_{m,i} \frac{(-1)^i} {i!}  \gamma^i (\gamma b+\delta)^{-2m-2n-i}  f_m(gb) X_{n,m+i}\\
=& \sum_{m=0}^\infty   \sum_{i=0} ^m \frac{(-1)^i}{i!} \gamma^i (\gamma b+\delta)^{-2m-2n+i}f_{m-i}(gb) X_{n,m}
\end{align*}
Therefore $L_n$ is fixed by $g$ if and only if
\begin{equation}  \label{5.21}
\sum_{i=0}^m \frac{(-1)^i}{i!} \gamma^i (\gamma b+\delta)^{-2m-2n+i} f_{m-i}(gb)=f_m(b), \;\;\;   \forall  m\geq 0,
\end{equation}  
which is equivalent to the following equations
\begin{equation}\label{5.22}
\sum_{i=0}^m \frac{1}{i!} \gamma^i (\gamma b+\delta)^{2m+2n-i} f_{m-i}(b)=f_m(gb), \;\;\;\forall m\geq 0.
\end{equation} 
Indeed, applying (\ref{5.21}) to the left hand side of (\ref{5.22}), we obtain 
\[
\sum_{i=0}^m \sum_{k=0}^{m-i }  \frac{(-1)^k}{i!k!} \gamma^{i+k} (\gamma b+\delta)^{i+k}f_{m-i-k}(gb)\\
=f_m(gb).\]
Hence (\ref{5.21}) implies (\ref{5.22}). Similarly we can show that (\ref{5.22}) also implies (\ref{5.21}).
Now let $\Phi(\tau ,X)=\sum_{m=0}^\infty f_m X^m$. Using (\ref{5.20}), we deduce that $\Phi$ is a Jacobi-like form if and only if (\ref{5.22}) holds.\qed

It is straight forward to show that the coefficients $\{c_{2n,m}f^{[m]}\}_{m}$ in (\ref{5.18}) satisfies (\ref{5.21}) thanks to (\ref{5.1}) and (\ref{5.13}). By Theorem \ref{theorem5.6},  for any $f\in \mathcal M_{2n}(\Gamma)$, the formal series
\[
\tilde f(\tau ,X):= \sum_{m=0}^\infty c_{2n,m} f^{[m]}X^m,
\]
satisfies the transformation property (\ref{5.20}). Especially when $f$ is a modular form,  $\tilde f$ is a Jacobi-like form. 
Let $h\in \mathcal M_{2k}(\Gamma)$. 
The coefficient of $X^s$ of the product
\[
\tilde f(\tau,-X) \tilde h (\tau,X)=\sum_{s=0}^\infty \sum_{m+l=s} c_{2n,m} c_{2k,l} f^{[m]} g^{[l]} X^s
\] 
is a meromorphic modular form of weight $2n+2k+2s$. Hence we define the generalized Rankin-Cohen bracket on meromorphic modular forms by
\begin{equation} \label{5.23}
[f,h]^\thicksim_{s}:=\sum_{m+l=s} c_{2n,m}c_{2k,l} f^{[m]} g^{[l]} \in \mathcal M_{2n+2k+2s}(\Gamma).
\end{equation}
\noindent Note that the generalized Rankin-Cohen brackets can be extended to meromorphic modular forms of odd weights by the same formula, and the brackets generalize the modified Rankin-Cohen brackets defined in \cite{D2}. Especially, when $f$ and $h$ are modular forms of positive weights, the generalized Rankin-Cohen bracket is exactly the usual Rankin-Cohen bracket.

\section{Acknowledgements}
The results reported here were obtained during the first author's stay at USTC during March to May 2022. 
The first author would like to thank Prof. Yongchang Zhu for discussions.
The research of the first author is supported by JSPS KAKENHI Grant Numbers 21H04993
and 23K19008.
The research of the second author is supported by National Natural Science Foundation of China No. 12171447.

\end{document}